\documentclass[letterpaper,10pt,twocolumn]{scrartcl}

\usepackage[top=18mm,bottom=23mm,left=18mm,right=18mm,columnsep=15pt]{geometry} 
\usepackage[format=plain,small,labelfont=bf,up,textfont=it,up]{caption} 
\usepackage{booktabs} 
\usepackage{type1cm}
\usepackage{lettrine} 
\usepackage{abstract} 
\usepackage{titlesec} 
\renewcommand\thesection{\Roman{section}} 
\renewcommand\thesubsection{\Alph{subsection}} 
\titleformat{\section}[block]{\scshape\centering}{\thesection.}{1em}{} 
\titleformat{\subsection}[block]{\it}{\thesubsection.}{1em}{} 

\usepackage[dvips]{graphicx}
\usepackage[usenames,dvipsnames]{xcolor}
\usepackage{tikz}
\usetikzlibrary{shapes,arrows,trees}
\usepackage{psfrag}
\usepackage{booktabs}

\usepackage[cmex10]{amsmath}
\usepackage{amssymb,amsfonts,amsthm}
\usepackage{mathrsfs}
\usepackage{dsfont}
\usepackage{cases}

\usepackage{algorithm}
\usepackage{algpseudocode}

 \usepackage[caption=false,font=footnotesize]{subfig}

\usepackage{fixltx2e}

\usepackage{url}

\usepackage[square,sort,numbers]{natbib}
\newtheorem{cor}{Corollary}
\newtheorem{lem}{Lemma}
\newtheorem{thm}{Theorem}

\newtheorem{rem}{Remark}

\definecolor{darkgreen}{rgb}{0,0.4,0}
\definecolor{darkgray}{rgb}{0.3,0.3,0.3}

\newcommand\TRANSP{{\top}}
\newcommand\pplus{+\!\!+}

\title{\fontsize{21pt}{16pt}\selectfont{\bf\normalfont A stabilizing iteration scheme for model predictive control\\[0.1cm] based on relaxed barrier functions}} 

\author{Christian~Feller~and~Christian~Ebenbauer %
\thanks{This work was supported by the Deutsche Forschungsgemeinschaft
(Emmy-Noether-Grant, Novel Ways in Control and Computation, EB 425/2-1, and Cluster of Excellence in Simulation Technology,
EXC 310/2).}\\[2mm] 
\normalsize Institute for Systems Theory and Automatic Control \\
\normalsize University of Stuttgart, Pfaffenwaldring 9, 70550 Stuttgart, Germany \\ 
\normalsize e-mail:\,\{christian.feller,ce\}@ist.uni-stuttgart.de 
\vspace{-5mm}
}
\date{}

\begin{document}
\twocolumn[
  \begin{@twocolumnfalse}
    \maketitle
\vspace*{-0.25cm}
\begin{abstract}
We propose and analyze a stabilizing iteration scheme for the algorithmic implementation of model predictive control for linear discrete-time systems. Polytopic input and state constraints are considered and handled by means of so-called relaxed logarithmic barrier functions.  
The required on-line optimization is based on warm starting and performs only a limited, possibly small, number of optimization algorithm iterations between two consecutive sampling instants. 
The optimization algorithm dynamics as well as the resulting suboptimality of the applied control input are taken into account explicitly in the stability analysis, and the origin of the resulting overall closed-loop system, consisting of state and optimization algorithm dynamics, is proven to be asymptotically stable. The corresponding constraint satisfaction properties are also analyzed.      
The theoretical results and a presented numerical example illustrate the fact that asymptotic stability as well as a satisfactory closed-loop performance can be achieved by performing only a single optimization algorithm iteration at each sampling step. 
\end{abstract}
\vspace*{0.6cm}
  \end{@twocolumnfalse}
]
\saythanks
%
%
\section{Introduction}
%
\lettrine[nindent=0em,lines=2]{M}odel predictive control (MPC)
is a control strategy that solves, at each sampling instant, a suitably chosen finite horizon open-loop optimal control problem in a receding horizon fashion. While there exist various theoretical results concerning stability properties of the closed-loop system, see, e.g.,~\cite{mayne00} and references therein, these are often based on the assumption that the exact optimal solution is available immediately when the next state measurement becomes available.  
However, this assumption might not be valid when considering applications in which the optimization algorithm operates on the same time scale as the system to be controlled, for example due to fast system dynamics and/or the use of low-cost hardware. 
In this case, special measures have to be taken in order to ensure stability properties, or even recursive feasibility, of the closed-loop system. \\
One possible approach is to speed up the optimization accordingly and then use again the aforementioned time-scale separation argument. In the case of linear systems, this might for example be achieved by making use of explicit MPC solutions~\cite{bemporadExplicitLQR,tondel03} or tailored, fast QP algorithms~\cite{rao98,ferreau08,wang10,koegel11,richter12,patrinos12,domahidi12,giselsson14}.
Another approach is to settle for approximate solutions and then take the suboptimality of the applied inputs explicitly into account within the stability analysis, see for example~\cite{scokaert99,mcgovern99,zeilinger11}, or to compute the optimal solution based on a one-step-ahead prediction~\cite{zavala09}. 
However, there are only few works which explicitly consider the \emph{dynamics} of the optimization algorithm and acknowledge the fact that the optimization typically needs to be performed in parallel to the evolving system dynamics. Interesting exceptions are given by the works~\cite{diehl05,diehl07}, in which the authors present, for the case of unconstrained nonlinear systems, a so-called real-time iteration scheme that performs only one Newton based optimizer update per sampling instant.
The authors rigorously analyze the combined dynamics consisting of system and optimizer state and provide a proof that the state of the corresponding overall closed-loop system asymptotically converges to zero over time.   
However, no input and state constraints are considered and the theoretical analysis is, as the authors themselves emphasize, carried out based on some technical assumptions that might be hard to verify for general practical applications.
Nevertheless, the results represent a theoretical underpinning of the proposed iteration scheme and nicely illustrate possible benefits that might be gained by combining ideas from the areas of MPC stability theory and optimization algorithms.     \\
In contrast to this, we focus in this paper on the linear MPC case with polytopic input and state constraints and show that, by making use of the relaxed barrier function based framework presented in~\cite{feller15a,feller15c}, an MPC iteration scheme can be derived that allows to guarantee desirable systems theoretic properties based on standard linear MPC assumptions that can be easily verified.
The iteration scheme is based on the idea to predict at each sampling instant the next system state, start the optimization procedure based on a suitable warm start solution, and then, at the end of the available computation time, simply apply the first element of the current, possibly suboptimal, optimizer state.  
Related to this scheme, the main contributions of this paper are as follows. 
First, we present an analysis of the resulting overall closed-loop system, consisting of plant and algorithm dynamics, in which we prove that, under standard assumptions, the origin of the combined state and optimizer dynamics will be asymptotically stable for any number of internal optimization algorithm iterations -- including the case of performing only a single optimizer update at each sampling step. 
Second, we also study the resulting constraint satisfaction properties and show that, for a suitable subset of feasible initial conditions, the input and state constraints can be satisfied with arbitrary tolerance -- including the case of exact constraint satisfaction -- if the underlying relaxation as well as the initialization of the iteration scheme are chosen in a suitable way.
Third, possible realizations of the underlying optimization algorithm are discussed based on gradient and Newton-type line search methods.
In summary, we thus present a barrier function based MPC iteration scheme for which we prove asymptotic stability as well as desirable constraint satisfaction properties of the closed-loop system, independently of both the actual realization of the underlying optimization algorithm and the number of performed optimization algorithm iterations.  
To the best of our knowledge, an MPC algorithm with such guaranteed systems theoretic and algorithmic properties has not been discussed so far in the literature.\\
Note that efficient optimization algorithms for (nonrelaxed) barrier function based linear MPC approaches have, for example, been discussed in~\cite{wright1997,rao98,wang10,domahidi12}. 
However, none of the aforementioned works analyzes the closed loop consisting of both system and optimization algorithm dynamics, and no proof of stability is given.
Nevertheless, some of the presented concepts may in the future be helpful for speeding up the required line search also in the relaxed barrier function based framework.  \\
The paper is organized as follows: necessary background on relaxed barrier function based MPC as well as some preliminary results are presented in Section~\ref{sec:prelimRes}; in Section~\ref{sec:mainRes}, we state and discuss our main results concerning the outlined MPC iteration scheme; in Section~\ref{sec:numExample} we illustrate some related aspects by means of a numerical example, before we state our conclusions in Section~\ref{sec:conclusion}. \\
Some remarks on the used notation are in place. $\mathbb{R}_+$, $\mathbb{R}_{+\!\!+}$, and $\mathbb{N}_+$ denote the sets of nonnegative real, strictly positive real, and strictly positive natural numbers. $\mathbb{S}^n_+$ and $\mathbb{S}^n_{\pplus}$ refer to the sets of positive semi-definite and positive definite matrices of dimension $n\in\mathbb{N}_+$. For any given matrix~$M$ or vector $v$, $M^i$ and $v^i$ refer to the $i$-th row or element, and $\lVert x \rVert_M:=(x^\TRANSP \! M x)^{1/2}$ for $M\in\mathbb{S}_+$. A polytope is defined as compact intersection of a finite number of halfspaces and $\mathcal{B}_r^n(z_0)$ is the $n$-dimensional Euclidean ball with radius $r$ and center $z_0$. For any set $S$, $S^\circ$ denotes the open interior. Moreover, $\mathds{1}:=[1 \ \cdots \ 1]^\TRANSP$.
A function $\alpha:\mathbb{R}_+ \to \mathbb{R}_+$ with $\alpha(0)=0$ is a $\mathcal{K}$-function if it is continuous and strictly increasing; it is a $\mathcal{K}_\infty$ function if $\alpha\in\mathcal{K}$ and in addition $\alpha(s)\to \infty$ for $s\to\infty$. A function $\beta:\mathbb{R}_+ \times \mathbb{R}_+ \to \mathbb{R}_+$ is a  $\mathcal{KL}$-function if $\beta(\cdot,t)\in \mathcal{K}$ for fixed $t\in\mathbb{R}_+$ while, for fixed $s\in\mathbb{R}_+$, the function $\beta(s,\cdot)$ is decreasing and $\beta(s,t)\to0$ as $t\to \infty$.

\section{Preliminary Results}\label{sec:prelimRes}
%
In the following, we introduce the problem setup as well as some necessary background on relaxed barrier function based MPC.  
\subsection{Problem setup}
We consider the control of linear time-invariant discrete-time systems of the form
\begin{equation}\label{eq:discreteSystem}
x(k+1)=Ax(k)+Bu(k)\, ,
\end{equation}
where $x(k) \in \mathbb{R}^n$ refers to the vector of system states and $u(k) \in \mathbb{R}^m$ refers to the vector of system inputs, both at time instant $k \geq 0$. Moreover, the matrices $A\in \mathbb{R}^{n\times n}$ and $B\in \mathbb{R}^{n\times m}$ describe the corresponding system dynamics, where we assume $(A,B)$ to be stabilizable. The control task is to regulate the system state to the origin while minimizing a given quadratic performance criterion and satisfying, whenever possible, polytopic state and input constraints of the form 
\begin{subequations} \label{eq:constraints}
\begin{align}
 x(k) \in \mathcal{X}&=\{x \in \mathbb{R}^n: C_{\mathrm{x}} x \leq d_{\mathrm{x}}\}\, , \\ 
 u(k) \in \mathcal{U}&=\{u \in \mathbb{R}^m: C_{\mathrm{u}} u \leq d_{\mathrm{u}}\} \, ,
\end{align} 
\end{subequations}
where $C_{\mathrm{x}}\in\mathbb{R}^{q_{\mathrm{x}}\times n}$, $C_{\mathrm{u}}\in\mathbb{R}^{q_{\mathrm{u}}\times m}$ and $d_{\mathrm{x}}\in\mathbb{R}^{q_{\mathrm{x}}}_{+\!\!+}$, $d_{\mathrm{u}}\in\mathbb{R}^{q_{\mathrm{u}}}_{+\!\!+}$ with $q_x, q_u\in\mathbb{N}_+$. Note that both $\mathcal{X}$ and $\mathcal{U}$ contain the origin in their interior. 
 In linear MPC, this problem setup is handled by solving, for current system state~$x=x(k)$ and finite prediction horizon~$N \in\mathbb{N}_+$, an open-loop optimal control problem of the form
\begin{subequations} \label{eq:OptProblemLin}
\begin{align}
J_N^{\ast}(x)&=\min_{\boldsymbol{u}} \sum_{k=0}^{N-1} \ell(x_k,u_k) + F(x_N) \\
\text{s.t.\ } \ &{x}_{k+1}=A{x}_{k}+B{u}_{k}, \ {x}_0=x \ , \\ 
&x_k \in \mathcal{X}, \ k=0,\dots, N-1, \ x_N \in \mathcal{X}_f \ ,\\
&u_k \in \mathcal{U}, \ k=0,\dots, N-1 \ ,
\end{align}
\end{subequations}
where the stage cost $\ell:\mathbb{R}^n\times \mathbb{R}^m\to \mathbb{R}_+$ and the terminal cost $F:\mathbb{R}^n\to\mathbb{R}_+$ are chosen as $\ell(x,u)=\lVert x\rVert_Q^2+\lVert u\rVert_R^2$ and $F(x)=\lVert x\rVert_P^2$ for appropriate weight matrices $Q\in \mathbb{S}^n_+$, $R \in \mathbb{S}^m_{\pplus}$, $P\in \mathbb{S}^n_{\pplus}$. Furthermore, $\boldsymbol{u}=\{u_0,\dots,u_{N-1}\}$ denotes the sequence of control inputs over the prediction horizon~$N$, while $\mathcal{X}_f$ refers to a closed and convex terminal constraint set that may be used to guarantee stability properties of the closed-loop system. Note that we make use of subindices to distinguish open-loop predictions~$x_k$,~$u_k$~from actual state and input trajectories~$x(k), u(k)$. The control law is obtained by solving~(\ref{eq:OptProblemLin}) at each sampling instant $k\geq 0$ and applying  $u(k)=u_0^*(x(k))$ in a receding horizon fashion, where $u_0^*$ is the first element of the optimal input sequence $\boldsymbol{u}^*$.
Different approaches on how to choose the corresponding parameters in order to guarantee stability properties of the closed-loop state dynamics are summarized in~\cite{mayne00} and references therein -- however, as discussed above, under the assumption that the exact optimal solution is available instantaneously.
%
\subsection{Relaxed Logarithmic Barrier Function Based MPC} \label{sec:relBarrierIntro}
The concept of relaxed logarithmic barrier function based MPC has been proposed recently in~\cite{feller15a}. 
Here, the term \emph{relaxed} relates to a logarithmic barrier function that is smoothly extended by a suitable penalizing term, which in general is assumed to be convex and globally defined~\cite{benTal92,nash94,hauser06}. For example, a relaxed logarithmic barrier function for the set $\mathbb{R}_+$ may be given by 
\begin{equation} \label{eq:relBarrierDef}
 \hat{B}(z)\!=\!\begin{cases}
             -\ln(z) &\! z>\delta \\
             \beta(z;\delta):=\frac{1}{2}\left[\left(\frac{z-2\delta}{\delta}\right)^2\!\!-\!1\right]\!-\ln(\delta) &\! z \leq \delta, 
            \end{cases}
\end{equation} 
where $\delta \in \mathbb{R}_{\pplus}$ denotes the so-called relaxation parameter and $\beta(\cdot;\delta)$ is a quadratic relaxing function~\cite{hauser06,feller15a}.
Obviously, the function $\hat{B}(\cdot)$ is globally defined, twice continuously differentiable, and convex. For a polytopic set of the form $\mathcal{P}=\{\xi \in\mathbb{R}^r : C \xi \leq d \}$ with $C\in\mathbb{R}^{q\times r}$, $d\in\mathbb{R}^q_{\pplus}$, we may simply use $\hat{B}_{\mathcal{P}}(\xi)=\sum_{i=1}^q \hat{B}(z_i(\xi))$ with $z_i(\xi)=-C^i \xi +d^i$ and $\hat{B}(\cdot)$ defined according to~(\ref{eq:relBarrierDef}). 
Positive definiteness of the barrier function~$\hat{B}_{\mathcal{P}}(\cdot)$ with respect to the origin can be achieved by making use of suitable recentering procedures, see~\cite{feller15b} and\cite{feller15a} for more details.
For example, when using a weighting based recentering approach, we obtain 
\begin{equation} \label{eq:relBarrierPol}
\hat{B}_{\mathcal{P}}(\xi)=\sum_{i=1}^q \big(1+w^i\big) \left(\hat{B}(-C^i \xi+d^i)+\ln(d^i) \right) \ ,
\end{equation}
which can easily be shown to be positive definite for any $0 < \delta \leq \min\{d^1,\dots,d^q\}$ and any $w \in \mathbb{R}^q_{+}$ satisfying $\nabla \hat{B}_{\mathcal{P}}(0)=C^\TRANSP \mathrm{diag}(\frac{1}{d^1},\dots,\frac{1}{d^q}) (\mathds{1}+w)=0$~\cite{feller15b}. 
For the rest of the paper we can simply think of the recentering as a modification which ensures positive definiteness of a (relaxed) barrier function while preserving its main characteristics.
The open-loop optimal control problem of relaxed barrier function based MPC then reads
\begin{subequations}\label{eq:OptProblemRelBarrier}
\begin{align}
\hat{J}_N^{\ast}(x)&=\min_{\boldsymbol{u}}\sum_{k=0}^{N-1} \hat{\ell}(x_k,u_k) + \hat{F}(x_N)   \\
\mbox{s.\,t.} \ \ & {x}_{k+1}=A{x}_{k}+B{u}_{k}, \ {x}_0=x\, .
\end{align}
\end{subequations}
Here, $\hat{\ell}(x,u):= \ell(x,u)+\varepsilon \hat{B}_{\mathrm{x}}(x) + \varepsilon \hat{B}_{\mathrm{u}}(u)$ is the modified stage cost, while $\hat{B}_{\mathrm{x}}(\cdot)$ and $\hat{B}_{\mathrm{u}}(\cdot)$ are recentered relaxed logarithmic barrier functions of the form~(\ref{eq:relBarrierPol}) for the polytopic constraint sets $\mathcal{X}$ and $\mathcal{U}$ in~(\ref{eq:constraints}) for a given relaxation parameter $\delta \in \mathbb{R}_{\pplus}$ and suitable weighting vectors $w_{\mathrm{x}}\in\mathbb{R}_+^{q_\mathrm{x}}$, $w_{\mathrm{u}}\in\mathbb{R}_+^{q_\mathrm{u}}$.
The positive scalar $\varepsilon \in \mathbb{R}_{+\!\!+}$ is the barrier function weighting parameter which determines the influence of the barrier function values on the overall cost function.
Note that the parameters  $\varepsilon$ and $\delta$ are assumed to take fixed values and are not decreased to zero, as it is for example done in conventional interior-point methods or the trajectory optimization approach discussed in~\cite{hauser06}.
Furthermore, $\hat{F}: \mathbb{R}^n \to \mathbb{R}_+$ denotes a suitable terminal cost function term which, as discussed in~\cite{feller15a}, is crucial for the stability properties of the closed-loop system. 
In the following, we consider a quadratic terminal cost function of the form
\begin{equation} \label{eq:quadraticF}
\hat{F}(x)=x^\TRANSP\!P\,x \, ,
\end{equation}
where $P\in\mathbb{S}^n_{\pplus}$ is a positive definite solution to the following algebraic Riccati equation
\begin{subequations} \label{eq:riccatiP}
\begin{align} 
K&=-\left(R+B^\TRANSP PB+\varepsilon M_{\mathrm{u}}\right)^{-1}B^\TRANSP PA  \\
P&=(A_{K}^{\TRANSP} P A_K^{}+K^\TRANSP \!(R+\varepsilon M_{\mathrm{u}})K+Q + \varepsilon M_{\mathrm{x}} \ ,
\end{align}
\end{subequations}
with $A_K:=A+BK$,
$M_{\mathrm{x}}:=\frac{1}{2\delta^2}C_{\mathrm{x}}^\TRANSP\mathrm{diag}\left(\mathds{1}+w_{\mathrm{x}}\right)C_{\mathrm{x}}^{}$ and $M_{\mathrm{u}}:=\frac{1}{2\delta^2}C_{\mathrm{u}}^\TRANSP\mathrm{diag}\left(\mathds{1}+w_{\mathrm{u}}\right)C_{\mathrm{u}}^{}$, respectively. For more details on this approach as well as on alternative choices for $\hat{F}(\cdot)$, we refer the reader to Section~IV.C of~\cite{feller15a}. 
By eliminating the linear equality constraints and using a vectorized representation, problem~(\ref{eq:OptProblemRelBarrier}) can be formulated as unconstrained minimization of a cost function of the form
\begin{equation}\label{eq:QPformulation}
\hat{J}_N(U,x)= \frac{1}{2}U^\TRANSP\! HU+x^{\!\TRANSP}\! FU + x^{\!\TRANSP}\!Yx + \varepsilon \hat{B}_{\mathrm{xu}}(U,x)\, ,
\end{equation}
where $U :=[u_0^\TRANSP \ \cdots \ u_{N-1}^\TRANSP]^\TRANSP \in \mathbb{R}^{Nm}$ is the stacked open-loop input vector and $\hat{B}_{\mathrm{xu}}: \mathbb{R}^{Nm} \times \mathbb{R}^n\to \mathbb{R}_+$ is a positive definite, convex, and continuously differentiable relaxed logarithmic barrier function for polytopic constraints of the form $GU\leq d+Ex$. 
The matrices $H \in \mathbb{S}^{Nm}_{\pplus}$, $F\in\mathbb{R}^{n\times Nm}$, $Y\in\mathbb{S}^{n}_+$, $G \in \mathbb{R}^{q \times Nm}$, $d\in \mathbb{R}^q$, and $E\in\mathbb{R}^{q \times n}$ with $q=N(q_{\mathrm{x}}+q_{\mathrm{u}})$ can be constructed from~(\ref{eq:OptProblemRelBarrier}) and the corresponding constraints by means of simple matrix operations, see Appendix~\ref{app:condensedForm}.
Due to the design of the relaxed barrier functions, the above cost function is twice continuously differentiable and globally defined. Furthermore, we can state the following results, for which a proof is given in Appendix~\ref{app:barrierProperties}.
\vspace*{-0.1cm}
\begin{lem}\label{lem:jProp}
 Let the cost function $\hat{J}_N:\mathbb{R}^{Nm}\times\mathbb{R}^n \to \mathbb{R}_+$ given in~(\ref{eq:QPformulation}) be associated to the relaxed barrier function based MPC approach based on~(\ref{eq:relBarrierDef})--(\ref{eq:riccatiP}). Then,
 \begin{itemize} \setlength{\itemsep}{2pt}
  \item[i)] $\exists\ \underline{\alpha}$, $\bar{\alpha} \in \mathcal{K}_\infty$ such that for all $(U,x) \in \mathbb{R}^{Nm}\times\mathbb{R}^n$ it holds that $\underline{\alpha}(\|(U,x)\|) \leq \hat{J}_N(U,x) \leq \bar{\alpha}(\|(U,x)\|)$;
  \item[ii)] $\hat{J}_N(U,x)$ is strongly convex in $U$, convex in $x$, and there exist $\sigma$, $L \in \mathbb{R}_{\pplus}$ such that for all $(U,x) \in \mathbb{R}^{Nm} \times \mathbb{R}^n$ it holds that
     \begin{subequations} \label{eq:Jprop}
   \begin{align}
    \sigma I \leq \nabla_U^2 \hat{J}_N(U,x) &\leq LI   \\[0.1cm]
    \hspace*{-0.47cm}\|\nabla_U\hat{J}_N(U_1,x)-\nabla_U\hat{J}_N(U_2,x)\| &\leq L\|U_1-U_2\|\, ; 
   \end{align}
  \end{subequations}
\item[iii)] $\exists\, k_f:\mathbb{R}^{Nm}\times\mathbb{R}^n\!\to\!\mathbb{R}^m$ with $k_f(0,0)=0$ such that for any $(U_0,x_0)$ with successor state $x_0^+\!:=\!Ax_0+B u_0$,
\begin{equation} \label{eq:jPropDecrease}
 \hat{J}_N(U_0^+,x_0^+)-\hat{J}_N(U_0,x_0) \leq - \hat{\ell}(x_0,u_0)\,, 
\end{equation}
where $\, U_0^+\!:=\!\begin{bmatrix}u_1^\TRANSP & \cdots & u_{N-1}^\TRANSP & k_f^\TRANSP(U_0,x_0) \end{bmatrix}^{\!\TRANSP}$\!.
 \end{itemize}
\end{lem}
\noindent Using property~$iii)$ of the above Lemma together with standard optimality arguments, it is easy to show that, for any initial condition $x(0)\in\mathbb{R}^n$ and any $k \in \mathbb{N}$, the feedback $u(k)=u_0^*(x(k))$ yields
\begin{equation}
 \hat{J}_N^*(x(k+1))-\hat{J}_N^*(x(k)) \leq -\hat{\ell}(x(k),\hat{u}_0^*(x(k))) \, ,
\end{equation} 
which implies that the origin of the closed-loop system will be globally asymptotically stable, cf. Theorem~5 in \cite{feller15a}. 
Thus, based on the above relaxed barrier function based formulation, an asymptotically stabilizing control input is given by
\begin{equation} \label{eq:piDef}
 u(k)=\Pi_0\, \hat{U}^*(x(k)),\ \ \  \Pi_0=\begin{bmatrix} I_m & 0 & \cdots & 0 \end{bmatrix}\, ,
\end{equation} 
where $\hat{U}^*(x)=\arg \min_U \hat{J}_N(U,x)$ is the state-dependent minimizer of~(\ref{eq:QPformulation}).
Note that Lemma~\ref{lem:jProp} implies that for any given $x\in\mathbb{R}^n$, the cost function $\hat{J}_N(\cdot,x)$ belongs to the class of strongly convex $C^2$ functions with Lipschitz continuous gradient, for which many well-known and efficient optimization algorithms exist, e.g., like fast and conjugated gradient methods, the heavy-ball, or the Newton method. 
For these algorithms, we may then in principle obtain bounds on the number of iterations, and hence on the computation time, that is required for a given desired level of suboptimality, cf.~\cite{mcgovern99,richter12}.  
However, in this paper we want to follow a different approach and propose an iteration scheme which allows not only to analyze the overall closed-loop system consisting of both system state and optimization algorithm dynamics in an integrated fashion but also to prove stability properties of this system for an arbitrary, and not necessarily a priori specified, number of optimization algorithm iterations -- including the case of performing only a single optimizer update at each sampling instant.
\section{Main Results}\label{sec:mainRes}
In this section we are going to present our main results, consisting of the outlined barrier function based MPC iteration scheme~(Section~\ref{subsec:iterationScheme}), an analysis concerning stability and constraint satisfaction properties of the corresponding closed-loop system~(Sections~\ref{subsec:stability} and \ref{subsec:constraints}), as well as a discussion on different realizations of the underlying iterative optimization algorithm~(Section~\ref{subsec:optOperator}).  
\subsection{The iteration scheme} \label{subsec:iterationScheme}
\begin{figure}[t]
\hspace{-0.1cm} \includegraphics[width=0.475\textwidth]{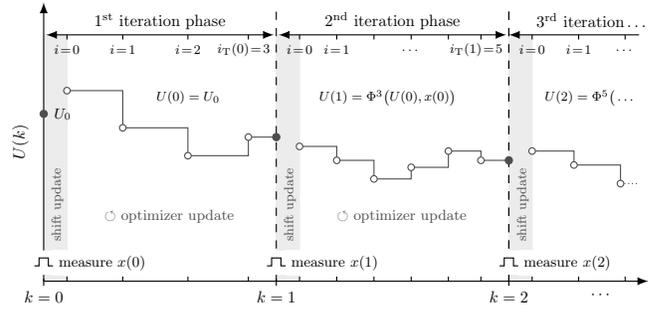}
 \caption{Schematic illustration of the proposed rbMPC iteration scheme with three optimizer iterations in the first and five optimizer iterations in the second iteration phase (scalar case). The resulting control law is $u(k)=\Pi_0\,U(k)$.}
 \label{fig:algoSchemeTime}
 \vspace*{-0.2cm}
\end{figure}

The rather intuitive idea underlying the following MPC iteration scheme is to predict at each sampling instant the next system state, start the optimization procedure for this predicted state based on a suitable warm start solution, and then, at the end of the available computation time, simply apply the first element of the current, possibly suboptimal, optimizer state.
A pseudo-code description of the iteration scheme is given in Algorithm~\ref{alg:rbMPCiteration} above, while a graphical illustration of the main idea can be found in Fig.~\ref{fig:algoSchemeTime}.
Aiming for a systems theoretic analysis, we may write the corresponding closed-loop system dynamics as
\begin{subequations} \label{eq:algoDyn}
 \begin{align}
  x(k+1)&=A\,x(k)+B \Pi_0\, U(k)\, , \hspace*{0.35cm} x(0)=x_0\, , \\
  U(k+1)&=\Phi^{i_{\mathrm{T}}(k)}\left(U(k),x(k)\right)\, , \hspace*{0.34cm} U(0)=U_0\, ,  
 \end{align}
\end{subequations}
where (\ref{eq:algoDyn}a) describes the dynamics of the controlled plant with $\Pi_0\in \mathbb{R}^{m\times Nm}$ being defined by~(\ref{eq:piDef}).
The dynamics of the controlled system state are coupled with the dynamics of the optimization algorithm, which are described by (\ref{eq:algoDyn}b). Here, $\Phi^{i_{\mathrm{T}}(k)}: \mathbb{R}^{Nm} \times \mathbb{R}^n \to \mathbb{R}^{Nm}$ denotes a suitable optimization algorithm operator that describes how the next optimizer state is computed at each sampling instant $k\in\mathbb{N}$ based on the current optimizer state $U(k)$ and the measured system state~$x(k)$ by performing $i_{\mathrm{T}}(k)\in\mathbb{N}$ optimization algorithm iterations. 
To be more precise, for given $i_{\mathrm{T}}(k)\in\mathbb{N}$, the operator $\Phi^{i_{\mathrm{T}}(k)}$ will be defined recursively as
\begin{subequations}\label{eq:recIteration}
\begin{align}
\Phi^{0}(U,x)&=\Psi_{\mathrm{s}}(U,x) \\
\Phi^{i}(U,x)&=\Psi_{\mathrm{o}}\!\left(\Phi^{i-1}(U,x),Ax+B\Pi_0U \right)\, ,
\end{align}
\end{subequations}
where a suitable warm start solution for the predicted next system state $x^+=,Ax+B\Pi_0U$ is  in the first step obtained via the \emph{shift operator} 
\begin{equation}\label{eq:shiftOperator}
\Psi_{\mathrm{s}}(U,x)=\begin{bmatrix} u_1^{\TRANSP} & u_2^{\TRANSP} & \cdots & u_{N-1}^{\TRANSP} & k^{\TRANSP} _f(U,x)\end{bmatrix}^{\TRANSP}\, ,
\end{equation} 
with $k_f:\mathbb{R}^{Nm}\times\mathbb{R}^n\!\to\!\mathbb{R}^m$ chosen according to Lemma~\ref{lem:jProp}. 
Furthermore, $\Psi_{\mathrm{o}}: \mathbb{R}^{Nm}\times \mathbb{R}^n \to \mathbb{R}^{Nm}$ denotes a so-called \emph{optimizer update operator}. We assume that
\begin{equation}\label{eq:optOperatorCond}
\hat{J}_N(\Psi_{\mathrm{o}}(U,x),x)-\hat{J}_N(U,x) \leq - \gamma(U,x)
\end{equation} 
for all $(U,x) \in \mathbb{R}^{Nm}\times\, \mathbb{R}^n$, where $\gamma: \mathbb{R}^{Nm} \times \mathbb{R}^n \to \mathbb{R}_{+}$ is a function satisfying $\gamma(U,x)\geq 0 \ \forall\, (U,x) \in \mathbb{R}^{Nm} \times \mathbb{R}^n$ and $\gamma(U,x)=0 \Leftrightarrow \nabla_U \hat{J}_N(U,x)=0$.
Different realizations of the operator~$\Psi_{\mathrm{o}}$ will be discussed later in Section~\ref{subsec:optOperator}, where we will also relate the function $\gamma(\cdot,\cdot)$ directly to the step size and the search direction of the respective optimization algorithm. For the moment, we may think of $\Phi^{i_{\mathrm{T}}(k)}$ as an operator that improves the current solution by applying a fixed number of suitable optimizer updates to the shifted warm start solution~(\ref{eq:shiftOperator}).
In the context of the rbMPC approach discussed above, we can exploit the fact that an explicit expression for a good warm start solution always exists, see property~$iii)$ of Lemma~\ref{lem:jProp}.
\begin{algorithm}[b]
\caption{rbMPC iteration scheme}\label{alg:rbMPCiteration}
\begin{algorithmic}[1]
\For{$k=0,1,2,\dots$}
\State obtain the current system state $x(k)$;
\State \textbf{if} $k=0$ \textbf{then} set $U(0)=U_0 \in \mathbb{R}^{Nm}$;
\State apply $u(k)=\Pi_0\, U(k)$ with $\Pi_0$ according to~(\ref{eq:piDef});
\State shift update: $\bar{U}_0(k)=\Psi_s\big(U(k),x(k)\big)$;
\For{$i=0,1,\dots,i_{\mathrm{T}}(k)-1$} optimizer update 
\State $\bar{U}_{i+1}(k)=\Psi_{\mathrm{o}}(\bar{U}_{i}(k),Ax(k)+B\Pi_0\,U(k))$;
\EndFor
\State set $U(k+1)=\bar{U}_{i_{\mathrm{T}}(k)}(k)$;
\EndFor
\end{algorithmic}
\end{algorithm}
In the following, we will show that, under rather mild assumptions, Algorithm~\ref{alg:rbMPCiteration} leads to global asymptotic stability of the overall closed-loop system~(\ref{eq:algoDyn}) and in addition allows us to give hard guarantees concerning the satisfaction of the now relaxed input and state constraints. 
\subsection{Stability properties of the closed-loop system} \label{subsec:stability}
Concerning the stability properties of the overall closed-loop system we can state the following results.
%
\begin{thm}\label{thm:stability}
Let the cost function $\hat{J}_N:\mathbb{R}^{Nm}\times\mathbb{R}^n \to \mathbb{R}_+$ given in~(\ref{eq:QPformulation}) be associated to the relaxed barrier function based MPC approach based on~(\ref{eq:relBarrierDef})--(\ref{eq:riccatiP}). 
Moreover, consider the closed-loop system~(\ref{eq:algoDyn}) with the operator $\Phi^{i_{\mathrm{T}}(k)}$ being defined according to~(\ref{eq:recIteration})--(\ref{eq:optOperatorCond}). Then, for any sequence $\boldsymbol{i}_{\mathrm{T}}=\{i_{\mathrm{T}}(0), i_{\mathrm{T}}(1), \dots\}$, the origin $(U,x)=(0,0)$ of system~(\ref{eq:algoDyn}) is globally asymptotically stable.
\end{thm}
\textit{\quad Proof: \ }
We know that $\underline{\alpha}(\|(U,x)\|) \leq \hat{J}_N(U,x) \leq \bar{\alpha}(\|(U,x)\|) \ \forall\, (U,x) \in \mathbb{R}^{Nm}\times\mathbb{R}^n$ for some $\underline{\alpha}$, $\bar{\alpha} \in \mathcal{K}_\infty$, see Lemma~\ref{lem:jProp}.
In the following, we show that $\hat{J}_N(U(k),x(k))$ will in addition always decrease to zero along trajectories of the closed-loop system. 
To this end, assume that we are given arbitrary initial conditions $(U_0,x_0) \in \mathbb{R}^{Nm} \times \mathbb{R}^n$. 
Consider now the next system state $x_0^+=Ax_0+B\Pi_0\,U_0$ as well as $U_0^+=\Psi_{\mathrm{s}}(U_0,x_0)$. With the shift operator $\Psi_{\mathrm{s}}(\cdot,\cdot)$ being defined by~(\ref{eq:shiftOperator}), it holds that 
\begin{equation}
 \hat{J}_N(U_0^+,x_0^+)-\hat{J}_N(U_0,x_0)\leq-\hat{\ell}(x_0,\Pi_0\,U_0)
\end{equation} 
due to property~$iii)$ of Lemma~\ref{lem:jProp}. Thus, $U_0^+$ can be seen as a warm start solution for the next system state that will lead to a decrease in the cost function.
Furthermore, condition~(\ref{eq:optOperatorCond}) ensures that the optimization algorithm update, i.e., applying the operator~$\Psi_{\mathrm{o}}(\cdot,\cdot)$, will lead to an even further decrease. In particular, $\hat{J}_N(\Phi^{i}(U_0,x_0),x_0^+)\leq \hat{J}_N(U_0^+,x_0^+)-\gamma(U_0^+,x_0^+) \leq \hat{J}_N(U_0,x_0)-\hat{\ell}(x_0,\Pi_0\,U_0)-\gamma(U_0^+,x_0^+)$ for any number $i\in \mathbb{N}_+$ of optimization algorithm iterations.  
As the above also holds for any $(U(k),x(k))$, this finally leads to 
\begin{align} \label{eq:stabProofDecrease}
  \hat{J}_N(U(k&+1),x(k+1))-\hat{J}_N(U(k),x(k))  \nonumber \\
 &\!\! \leq -\hat{\ell}(x(k),\Pi_0\,U(k)) -\gamma(U^+(k),x^+(k))   
\end{align} 
for all $k\in\mathbb{N}$, where $U^+(k):=\Psi_{\mathrm{s}}(U(k),x(k))$ and $x^+(k)=Ax(k)+B\Pi_0\, U(k)$. As both $\hat{\ell}(\cdot,\cdot)$ and $\gamma(\cdot,\cdot)$ are nonnegative, this implies that $(U,x)=0$ is stable. In particular $\|(U(k),x(k))\|\leq \underline{\alpha}^{-1}( \bar{\alpha}(\|(U(0),x(0))\|))$ by property $i)$ of Lemma~\ref{lem:jProp}.
Furthermore, the only solution that can stay identically in
$\Omega=\{(U,x)\in \mathbb{R}^{Nm}\times\mathbb{R}^n: \hat{\ell}(x,\Pi_0\,U)+\gamma(U^+,x^+) =0 \}$ 
is given by $(U,x)=0$, by which global asymptotic stability of the origin follows from the Barbashin-Krasovskii Theorem~\cite{halanay2000}.
The fact that $\Omega=\{0\}$ can be shown with the following arguments. First of all, positive definiteness of $\hat{\ell}(\cdot,\cdot)$ implies $x=0$ as well as $\Pi_0\, U=0$. As a direct consequence, $x^+=Ax+B\Pi_0\,U=0$. Moreover, $\gamma(U^+,x^+)=0$ implies that $U^+=\hat{U}^*(x^+)$ and, hence, $U^+=\hat{U}^*(0)=0$. As $U^+=\Psi_{\mathrm{s}}(U,x)$ with $\Psi_{\mathrm{s}}(\cdot,\cdot)$ defined according to (\ref{eq:shiftOperator}), it follows that $u_i=0$ for $i=1,\dots,N-1$. Together with $\Pi_0\,U=u_0=0$, this finally shows that $U=0$. \hfill $\blacksquare$ \\[0.2cm]
%
A direct consequence of Theorem~\ref{thm:stability} is that a fixed number of optimization algorithm updates may be chosen by the user without jeopardizing the stability properties of the overall closed-loop system -- including the case of performing only one optimizing step at each sampling instant. On the other hand, $i_{\mathrm{T}}(k)$ may also be time-varying and implicitly defined as the number of iterations after which the optimization algorithm is forced to return a solution. 
In this light, Algorithm~\ref{alg:rbMPCiteration} can be seen as a stabilizing \emph{anytime MPC algorithm} in the sense that it returns a valid and asymptotically stabilizing control input after each of its internal iterations, cf.~\cite{bemporad15}. 
Note, however, that no saturation or truncation of input or state variables has been considered, cf. Remark 4 in~\cite{feller15b}. 
\begin{rem}
We want to emphasize that procedures like computing the control input with a one-step-ahead prediction or performing only a limited number of optimization algorithm iterations are well-known in the literature, see~\cite{zavala09,scokaert99,wang10} and, in particular,~\cite{diehl05,diehl07}. To the best of the authors knowledge, however, no iteration scheme of the form~(\ref{eq:algoDyn}) with stability properties as in Theorem~\ref{thm:stability} has been presented so far. 
Note also that the presented stability results follow directly from the properties of the relaxed barrier function based MPC approach and that no additional assumptions are needed.
\end{rem}
\begin{rem}
 Note that Theorem~\ref{thm:stability} implies that if the optimizer dynamics are initialized with
$\|U_0\|\leq \alpha_0(\|x_0\|)$ for some $\alpha_0 \in \mathcal{K}_\infty$, then the controlled system state dynamics will be asymptotically stable in the sense that there exists $\beta_{\mathrm{x}}\in\mathcal{KL}$ such that 
$\|x(k)\|\leq \beta_{\mathrm{x}}(\|x_0\|,k)$ for all $x(0)=x_0\in\mathbb{R}^n$. \! \!\footnote{If $\beta \in \mathcal{KL}$ is the $\mathcal{KL}$ function related to the stable overall system dynamics~(\ref{eq:algoDyn}), then $\beta_{\mathrm{x}}(r,s)= \beta((\alpha_0 + \mathrm{id})(r),s)$.}
Suitable initializations with this property are for example given by $U_0=\hat{U}^*(x_0)$, i.e., the optimal relaxed barrier function based input for $x_0$, by $U_0=\bar{K} x_0$ with $\bar{K}\in\mathbb{R}^{Nm\times n}$, or simply by~$U_0=0$. 
\end{rem}
%
%
\subsection{Guarantees on closed-loop constrain satisfaction} \label{subsec:constraints}
When applying the above MPC iteration scheme, the input and state constraints in~(\ref{eq:constraints}) may be violated in closed-loop operation due to the relaxation of the underlying barrier functions. However, depending on the choice of the relaxation parameter~$\delta$ and the initialization of the iteration scheme, we can nevertheless state several interesting results concerning the constraint satisfaction properties of the closed-loop system. 
We begin our studies with the following Theorem, which provides an upper bound on the maximal constraint violation that is, in addition, shown to be monotonically decreasing over time.  
\begin{thm} \label{thm:maxConstrViolation}
Let the cost function $\hat{J}_N:\mathbb{R}^{Nm}\times\mathbb{R}^n \to \mathbb{R}_+$ given in~(\ref{eq:QPformulation}) be associated to the relaxed barrier function based MPC approach based on~(\ref{eq:relBarrierDef})--(\ref{eq:riccatiP}). 
Consider system~(\ref{eq:algoDyn}) with the operator $\Phi^{i_{\mathrm{T}}(k)}$ being defined according to~(\ref{eq:recIteration})--(\ref{eq:optOperatorCond}).
Moreover, let $\{x(0),x(1),\dots \}$ and $\{u(0),u(1),\dots \}$ with $u(k)=\Pi_0\,U(k)$ denote the resulting closed-loop state and input trajectories. 
Then, for any initialization $(U_0,x_0) \in \mathbb{R}^{Nm} \times \mathbb{R}^n$, any sequence $\boldsymbol{i}_{\mathrm{T}}=\{i_{\mathrm{T}}(0), i_{\mathrm{T}}(1), \dots\}$, and any $k\in\mathbb{N}$ it holds that 
\begin{equation} \label{eq:maxConstrViolThm}
    C_{\mathrm{x}} x(k)\leq d_{\mathrm{x}} + \hat{z}_{\mathrm{x}}(k), \quad  C_{\mathrm{u}} u(k)\leq d_{\mathrm{u}} + \hat{z}_{\mathrm{u}}(k),
\end{equation} 
where the elements of the maximal constraint violation vectors $\hat{z}_{\mathrm{x}}(k) \in \mathbb{R}^{q_{\mathrm{x}}}$ and $\hat{z}_{\mathrm{u}}(k) \in \mathbb{R}^{q_{\mathrm{u}}}$ decrease monotonically over time and are given by~(\ref{eq:maxConstrViol}). Furthermore, there exists a finite $k_0\in \mathbb{N}$ such that $\hat{z}_{\mathrm{x}}(k)\leq 0$, $\hat{z}_{\mathrm{u}}(k)\leq 0$ $\forall k\geq k_0$.
\end{thm}
\textit{\quad Proof: \ }
The proof consists of two parts. First we show that for any $k\in\mathbb{N}$, the relaxed barrier functions satisfy 
\begin{subequations}\label{eq:BxBuBounds}
\begin{align}
 \varepsilon\hat{B}_{\mathrm{x}}(x(k))&\leq \hat{J}_N(U(k),x(k))-x(k)^{\!\TRANSP}\! P_{\mathrm{uc}}^* x(k)\, , \\
 \varepsilon\hat{B}_{\mathrm{u}}(u(k))&\leq \hat{J}_N(U(k),x(k))-x(k)^{\!\TRANSP}\! P_{\mathrm{uc}}^* x(k)\, ,
\end{align}
\end{subequations}
where $P_{\mathrm{uc}}^* \in \mathbb{S}^n_{\pplus}$ is the solution to the discrete-time algebraic Riccati equation~(DARE) related to the infinite-horizon LQR problem for the setup $(A,B,Q,R)$. In a second step we then derive an upper bound for the maximal constraint violations and use the decrease of the value function to show the claimed monotonicity result. \\
\emph{Part 1.} 
Based on (\ref{eq:stabProofDecrease}) and positive definiteness of $\gamma(\cdot,\cdot)$, we know that
 $\hat{J}_N(U(i+1),x(i+1))-\hat{J}_N(U(i),x(i)) \leq -\hat{\ell}(x(i),u(i))$ for any $i\in\mathbb{N}$ and any initialization $(U_0,x_0) \in \mathbb{R}^{Nm} \times \mathbb{R}^n$, which in particular implied that $\lim_{i\to \infty}\hat{J}_N(U(i),x(i))=0$. Summing up over all future sampling instants beginning at $i=k$ and using a telescoping sum on the left hand side, we get $\hat{J}_N(U(k),x(k)) \geq \sum_{i=k}^\infty \hat{\ell}(x(i),u(i))$ for any $k \in \mathbb{N}$. Furthermore, it holds that $\sum_{i=k}^\infty \hat{\ell}(x(i),u(i))$ $\geq x(k)^\TRANSP \! P_{\mathrm{uc}}^*\, x(k) +\varepsilon \sum_{i=k}^\infty \hat{B}_{\mathrm{x}}(x(i))+\hat{B}_{\mathrm{u}}(u(i))$ due to the definition of $\hat{\ell}(\cdot,\cdot)$ and the optimality of the unconstrained infinite-horizon LQR solution. In combination, this yields
 $ \hat{J}_N(U(k),x(k)) \geq x(k)^\TRANSP \! P_{\mathrm{uc}}^* \,x(k)+\varepsilon \sum_{i=k}^\infty \hat{B}_{\mathrm{x}}(x(i))+\hat{B}_{\mathrm{u}}(u(i))$
 and finally, since all terms in the sum on the right hand side are positive definite, 
 $\varepsilon \hat{B}_{\mathrm{x}}(x(k)) \leq   \hat{J}_N(U(k),x(k))- x(k)^\TRANSP \! P_{\mathrm{uc}}^*\, x(k)$ as well as 
 $ \varepsilon \hat{B}_{\mathrm{u}}(u(k)) \leq   \hat{J}_N(U(k),x(k))- x(k)^\TRANSP \! P_{\mathrm{uc}}^*\, x(k)$ $\forall\,k \in \mathbb{N}$.\\
\emph{Part 2.} For ease of notation, let us define \vspace*{-0.15cm}
\begin{equation} \label{eq:alphaEq}
\hat{\alpha}(k):= \hat{J}_N(U(k),x(k))-x(k)^\TRANSP P_{\mathrm{uc}}^*\, x(k)^{} \, .
\end{equation}
From~(\ref{eq:BxBuBounds}) and~(\ref{eq:alphaEq}) it follows immediately that upper bounds for the maximal violations of state and input constraints are given by
\begin{subequations}\label{eq:maxConstrViol}
\begin{align}
 \hat{z}_{x}^i(k)= &\max_{\xi} \big\lbrace C_{\mathrm{x}}^i \xi -d_{\mathrm{x}}^i\, \big\vert \ \varepsilon \hat{B}_{\mathrm{x}}(\xi)\leq \hat{\alpha}(k) \big\rbrace\, ,  \\
  \hat{z}_{u}^j(k)= &\max_{v} \big\lbrace C_{\mathrm{u}}^j v -d_{\mathrm{u}}^j\, \big\vert  \ \varepsilon \hat{B}_{\mathrm{u}}(v)\leq \hat{\alpha}(k)\big\rbrace\, , 
\end{align}
\end{subequations}
where $i=1,\dots,q_{\mathrm{x}}$ and $j=1,\dots,q_{\mathrm{u}}$
 This shows that~(\ref{eq:maxConstrViolThm}) holds with the elements of $\hat{z}_{\mathrm{x}}(k) \in \mathbb{R}^{q_{\mathrm{x}}}$ and $\hat{z}_{\mathrm{u}}(k) \in \mathbb{R}^{q_{\mathrm{u}}}$ given by~(\ref{eq:maxConstrViol}). 
 Note that the optimization problems in~(\ref{eq:maxConstrViol}) are convex as the involved barrier functions are convex. 
 Moreover, as can be seen from the following arguments, $\hat{\alpha}(k)$ is monotonically decreasing and $\lim_{k\to\infty} \hat{\alpha}(k)=0$. On the one hand, the decrease in the cost function implies that $\hat{\alpha}(k+1)-\hat{\alpha}(k)=\hat{J}_N(U(k+1),x(k+1))-\hat{J}_N(U(k),x(k))-x(k+1)^\TRANSP\! P_{\mathrm{uc}}^*\, x(k+1)^{}+x(k)^\TRANSP\! P_{\mathrm{uc}}^*\, x(k)^{}$ $\leq -\hat{\ell}(x(k),u(k))-x(k+1)^\TRANSP\! P_{\mathrm{uc}}^*\, x(k+1)^{}+ x(k)^\TRANSP\! P_{\mathrm{uc}}^*\, x(k)^{}$. On the other hand, $\hat{\ell}(x(k),u(k))= \ell(x(k),u(k))+\varepsilon \hat{B}_{\mathrm{x}}(x(k))+\varepsilon\hat{B}_{\mathrm{u}}(u(k))$, while the principle of optimality ensures that ${\ell}(x(k),u(k))+x(k+1)^\TRANSP\! P_{\mathrm{uc}}^*\, x(k+1)^{} \geq x(k)^\TRANSP\! P_{\mathrm{uc}}^*\, x(k)^{}$. Thus, 
\begin{equation} \label{eq:alphaDecrease}
 \hat{\alpha}(k+1)-\hat{\alpha}(k) \leq -\varepsilon (\hat{B}_{\mathrm{x}}(x(k))+\hat{B}_{\mathrm{u}}(u(k)))\, ,
 \end{equation} 
which proves our claims since $\hat{B}_{\mathrm{x}}(\cdot)$ and $\hat{B}_{\mathrm{u}}(\cdot)$ are positive definite.
By these properties of $\hat{\alpha}(\cdot)$, it follows immediately that $\hat{z}^i_{\mathrm{x}}(k)$ and $\hat{z}^j_{\mathrm{u}}(k)$ given by~(\ref{eq:maxConstrViol}) also decrease monotonically over time. Finally, as both $x(k)$ and $u(k)$ converge asymptotically to the origin~(see Theorem~\ref{thm:stability}), there exists a finite $k_0\in \mathbb{N}$ such that $\hat{z}_{\mathrm{x}}(k)\leq 0$, $\hat{z}_{\mathrm{u}}(k)\leq 0$ $\forall\, k\geq k_0$. \hfill $\blacksquare$ 
\begin{rem}
 Note that the above results imply that the maximal violations of the state and input constraints will always be upper bounded by $\hat{z}_{\mathrm{x}}(0)$ and $\hat{z}_{\mathrm{u}}(0)$, respectively, which can for a given initialization $(U_0,x_0)$ be computed by solving the convex optimization problems stated in~(\ref{eq:maxConstrViol}) for $\hat{\alpha}(0):= \hat{J}_N(U_0,x_0)-x_0^\TRANSP P_{\mathrm{uc}}^*\, x_0^{}$.
\end{rem}
Thus, for any choice of the relaxation parameter $\delta$, an upper bound for the maximal possible constraint violation can be computed a priori. What is more, the following result, in fact, shows that we can in many cases ensure approximate (or even exact) constraint satisfaction by choosing relaxation parameter and initialization accordingly. In particular, the core message of the following theorem is, that at least for a suitable set of feasible initial conditions, the input and state constraints can be satisfied with any desired tolerance if we choose $\delta$ small enough and initialize the iteration scheme good enough, i.e., close enough to the optimal solution.  \vspace*{-0.1cm}
\begin{thm} \label{thm:xNdelta0}
Let the cost function $\hat{J}_N:\mathbb{R}^{Nm}\times\mathbb{R}^n \to \mathbb{R}_+$ given in~(\ref{eq:QPformulation}) be associated to the relaxed barrier function based MPC approach based on~(\ref{eq:relBarrierDef})--(\ref{eq:riccatiP}). 
Consider system~(\ref{eq:algoDyn}) with the operator $\Phi^{i_{\mathrm{T}}(k)}$ being defined according to~(\ref{eq:recIteration})--(\ref{eq:optOperatorCond}).
Let $\{x(0),x(1),\dots \}$ and $\{u(0),u(1),\dots \}$ with $u(k)=\Pi_0\,U(k)$ denote the resulting closed-loop state and input trajectories. 
Furthermore, assume $(A,B)$ to be controllable and let the set $\mathcal{X}_N$ be defined as 
\begin{align}\label{eq:xN_def_quadr}
  \mathcal{X}_N = \{x \in \mathbb{R}^n : \exists\, U\in\mathbb{R}^{Nm},  u_k\in \mathcal{U}, \text{\ such that\ } \\ 
  x_k(U,x) \in \mathcal{X}, x_N(U,x) = 0\}\, , \nonumber
\end{align}
where $x_0=x$ while $x_k(U,x)$ for $k=1,\dots,N$ are given according to~(\ref{eq:xkDyn}).
Then, for any compact set of initial conditions $\mathcal{X}_0 \subseteq \mathcal{X}_N^\circ$ and any $\hat{z}_{\mathrm{x},\text{tol}} \in \mathbb{R}^{q_x}_+$, $\hat{z}_{\mathrm{u},\text{tol}} \in \mathbb{R}^{q_u}_+$, there exists $\bar{\delta}_0 \in \mathbb{R}_{\pplus}$ and $\mu \in \mathcal{K}_\infty$ such that for any relaxation parameter $0<\delta \leq \bar{\delta}_0$, any $x_0 \in \mathcal{X}_0$, any initialization $U_0 \in \mathcal{B}_{\mu(\delta)}^{Nm}(\hat{U}^*(x_0))$, any $\boldsymbol{i}_{\mathrm{T}}=\{i_{\mathrm{T}}(0), i_{\mathrm{T}}(1), \dots\}$, and any $k\in\mathbb{N}$, it holds that 
\begin{equation} \label{eq:maxConstrViolTolerance}
    C_{\mathrm{x}} x(k)\leq d_{\mathrm{x}} + \hat{z}_{\mathrm{x,tol}}\,, \quad  C_{\mathrm{u}} u(k)\leq d_{\mathrm{u}} + \hat{z}_{\mathrm{u,tol}}\, .
\end{equation}
\end{thm}
\textit{\ Proof: }
We proof the result for the case of no constraint violation, i.e. $\hat{z}_{\mathrm{x,tol}}=0$, $\hat{z}_{\mathrm{u,tol}}=0$, which comprises all cases with positive tolerances. The proof consists of two parts.
First, we show that there exists for any $\delta \in \mathbb{R}_{\pplus}$ a compact and nonempty set $\hat{\mathcal{Z}}_N(\delta)\!\subseteq \mathbb{R}^{Nm} \times \mathbb{R}^n$ such that for any $(U_0,x_0) \in \hat{\mathcal{Z}}_N(\delta)$, the state and input constraints will be satisfied also for all future steps, i.e., (\ref{eq:maxConstrViolTolerance}) will hold with $\hat{z}_{\mathrm{x,tol}}=0$, $\hat{z}_{\mathrm{u,tol}}=0$.
Then, we show in the second part that for any compact set $\mathcal{X}_0 \subseteq \mathcal{X}_N^\circ$ there exists  $\bar{\delta}_0\in\mathbb{R}_{\pplus}$ and $\mu \in \mathcal{K}_\infty$ such that if $\delta\leq \bar{\delta}_0$, $x_0 \in \mathcal{X}_0$, and $U_0 \in \mathcal{B}_{\mu(\delta)}^{Nm}(\hat{U}^*(x_0))$, then $(U_0,x_0) \in \hat{\mathcal{Z}}_N(\delta)$. \\[0.1cm]
\emph{Part 1.} 
Let us for a given relaxation parameter $\delta \in \mathbb{R}_{\pplus}$ define the scalars $\bar{\beta}_{\mathrm{x}}(\delta), \bar{\beta}_{\mathrm{u}}(\delta), \bar{\beta}(\delta) \in \mathbb{R}_{\pplus}$ as
\begin{subequations} \label{eq:betaXUEq}
\begin{align}
\bar{\beta}_{\mathrm{x}}(\delta)&:=\min_{i,\xi}\{\hat{B}_{\mathrm{x}}(\xi)\vert \, {C_{\mathrm{x}}^i \xi=d_{\mathrm{x}}^i}\}, \ \, i=1,\dots,q_{\mathrm{x}} \, , \\
\bar{\beta}_{\mathrm{u}}(\delta)&:=\min_{j,v}\{\hat{B}_{\mathrm{u}}(v) \vert \, {C_{\mathrm{u}}^j v=d_{\mathrm{u}}^j}\}, \ \, j=1,\dots,q_{\mathrm{u}} \, , \\
\bar{\beta}(\delta)&=\min\{\bar{\beta}_{\mathrm{x}}(\delta),\bar{\beta}_{\mathrm{u}}(\delta)\} \, .
\end{align}
\end{subequations}
Note that $\bar{\beta}_{\mathrm{x}}(\delta)$ and $\bar{\beta}_{\mathrm{u}}(\delta)$ can be interpreted as lower bounds for the values that are attained by the relaxed barrier functions on the boundaries of the constraint sets $\mathcal{X}$ and $\mathcal{U}$, respectively.  
As a consequence, the $\bar{\beta}(\delta)$-sublevel sets of $\hat{B}_{\mathrm{x}}(\cdot)$ and $\hat{B}_{\mathrm{u}}(\cdot)$ will always be contained within the sets $\mathcal{X}$, $\mathcal{U}$ -- see Lemma~\ref{lem:relBarrierLevelSets} in the Appendix.
Based on this observation and inspired by the proof of Theorem~\ref{thm:maxConstrViolation}, we introduce
\begin{equation}\label{eq:alphaUX}
\hat{\alpha}(U,x):= \hat{J}_N(U,x)-x^\TRANSP P_{\mathrm{uc}}^*\, x^{}
\end{equation}
and define the set $\hat{\mathcal{Z}}_{N}(\delta) \subset \mathbb{R}^{Nm} \times \mathbb{R}^n$ as
 \begin{equation}
  \hat{\mathcal{Z}}_{N}(\delta)\!:=\!\left\lbrace (U,x) \in \mathbb{R}^{Nm}\!\times\!\mathbb{R}^n\, |\, \hat{\alpha}(U,x) \leq \varepsilon \bar{\beta}(\delta) \right\rbrace,
 \end{equation} 
 where $\bar{\beta}(\delta)$ is given by~(\ref{eq:betaXUEq}) and $P_{\mathrm{uc}}^* \in \mathbb{S}^n_{\pplus}$ is again the solution to the DARE related to the infinite-horizon LQR problem for $(A,B,Q,R)$. Suppose now that $(U_0,x_0) \in \hat{\mathcal{Z}}_{N}(\delta)$. Based on the proof of Theorem~\ref{thm:maxConstrViolation}, we know that $\varepsilon \hat{B}_{\mathrm{x}}(x(k))\leq \hat{\alpha}(U(k),x(k))$, $\varepsilon \hat{B}_{\mathrm{u}}(u(k))\leq \hat{\alpha}(U(k),x(k))$, see~(\ref{eq:BxBuBounds}). In addition, $\hat{\alpha}(U(k),x(k)) \leq \hat{\alpha}(U_0,x_0)$ for all $k \in \mathbb{N}$, cf.~(\ref{eq:alphaDecrease}). So if $(U_0,x_0) \in \hat{\mathcal{Z}}_{N}(\delta)$, then
 $\hat{B}_{\mathrm{x}}(x(k))\leq \bar{\beta}(\delta)$, $\hat{B}_{\mathrm{u}}(u(k))\leq \bar{\beta}(\delta)$, which by the definition of $\bar{\beta}(\delta)$ and Lemma~\ref{lem:relBarrierLevelSets} implies that $x(k) \in \mathcal{X}$, $u(k)\in\mathcal{U}$ for all $k\in\mathbb{N}$. Thus, for given $\delta\in\mathbb{R}_{\pplus}$,~(\ref{eq:maxConstrViolTolerance}) holds with $\hat{z}_{\mathrm{x,tol}}=0$, $\hat{z}_{\mathrm{u,tol}}=0$ whenever $(U_0,x_0) \in \hat{\mathcal{Z}}_{N}(\delta)$. \\[0.1cm]
\emph{Part 2.} 
In the following, we are going to prove our second claim, i.e., the existence of $\bar{\delta}_0 \in \mathbb{R}_{\pplus}$ and $\mu \in \mathcal{K}_\infty$ with the properties specified above.
As we need to vary the relaxation parameter, we make the influence of $\delta$ on the cost function more explicit.
In particular, we will use $\hat{J}_N(U,x;\delta)$ to denote the value of the cost function for $(U,x)\in\mathbb{R}^{Nm}\times\mathbb{R}^n$ and a given relaxation parameter $\delta\in \mathbb{R}_{\pplus}$. In accordance with~(\ref{eq:alphaUX}), we further define $\hat{\alpha}(U,x;\delta):=\hat{J}_N(U,x;\delta)-x^\TRANSP P_{\mathrm{uc}}^*\, x^{}$. Based on this notation, it obviously holds that
\begin{equation} \label{eq:zNdeltaEquivalence}
 (U_0,x_0)\in \hat{\mathcal{Z}}_N(\delta)\, \Leftrightarrow \, \hat{\alpha}(U_0,x_0;\delta) \leq \varepsilon \bar{\beta}(\delta)\, .
\end{equation} 
 We now proceed as follows: first we show in \emph{a)} that for a \emph{fixed} $x_0 \in \mathcal{X}_0 \subseteq \mathcal{X}_N^\circ$ and a corresponding suitably chosen $U_0$ there exists a $\delta_0'(x_0)$ such that~(\ref{eq:zNdeltaEquivalence}) holds for all $\delta \leq \delta_0'(x_0)$;  in part \emph{b)}, we then show the claimed result for $\mathcal{X}_0$ by \emph{uniformly} choosing $\bar{\delta}_0=\min_{x_0\in\mathcal{X}_0}{\delta}'_0(x_0)$.\\
 \emph{a)\,} Assume $x_0\in \mathcal{X}_0\subseteq \mathcal{X}_N^\circ$ and let $\tilde{U}^*(x_0)$ and $\tilde{J}_N^*(x_0)=\tilde{J}_N(\tilde{U}^*(x_0),x_0)$ denote the optimal input vector and value function for a \emph{nonrelaxed} barrier function based problem formulation with \emph{additional} constraint $x_N=0$, see Lemma~\ref{lem:continuityRes} in the Appendix. 
  Note that the additional constraint $x_N=0$ is required in order to ensure that $\tilde{J}_N^*(x_0)$ does not depend on~$\delta$ via the quadratic terminal cost defined in~(\ref{eq:quadraticF}) and (\ref{eq:riccatiP}).
 By definition, $\tilde{U}^*(x_0)$ exists for any $x_0\in \mathcal{X}_0$ and results in strictly feasible state and input sequences $\{\tilde{x}_0,\tilde{x}_1,\dots \tilde{x}_N\}$ and $\{\tilde{u}_0,\tilde{u}_1,\dots \tilde{u}_{N-1}\}$, respectively. 
 Now define $\delta_0(x_0):=\min_{i,j,k}\{-C_{\mathrm{x}}^i \tilde{x}_k+d_{\mathrm{x}}^i$, $-C_{\mathrm{u}}^j \tilde{u}_k+d_{\mathrm{u}}^j\}$, where $i = 1,\dots,q_{\mathrm{x}}$, $j=1,\dots,q_{\mathrm{u}}$, $k=0,\dots,N-1$, which characterizes the minimal distance of the nonrelaxed open-loop trajectories to the boundaries of the respective constraint sets.
Due to (\ref{eq:relBarrierDef}), the cost function values of relaxed and nonrelaxed formulation will be identical, i.e., $\hat{J}_N(\tilde{U}^*(x_0),x_0;\delta)=\tilde{J}_N^*(x_0)$, if we choose $\delta \leq \delta_0(x_0)$, see also~\cite{hauser06}. In addition, it always holds by definition that $\hat{J}_N^*(x_0;\delta) \leq \hat{J}_N(\tilde{U}^*(x_0),x_0;\delta)$ and, hence, $\hat{J}_N^*(x_0;\delta) \leq \tilde{J}_N^*(x_0)$ for all $\delta \leq \delta_0(x_0)$.
 As the smooth and convex problem formulation ensures that $\tilde{U}^*(x_0)$ will be continuous, see Lemma~\ref{lem:continuityRes} in the Appendix, $\delta_0(x_0)$ will be a continuous function of $x_0$. \\
 Next, we derive an upper bound for the cost function $\hat{J}_N(U_0,x_0;\delta)$ for input vectors $U_0$ in the neighborhood of~$\hat{U}^*(x_0)$. 
 Assume that $\delta \leq \delta_0(x_0)$ and consider $U_0 \in \mathcal{B}_{\mu_0}^{Nm}(\hat{U}^*(x_0))$ for arbitrary $\mu_0\in\mathbb{R}_{\pplus}$, i.e. $U_0=\hat{U}^*(x_0)+\mu_0 d$, $\|d\|\leq 1$. Due to Taylor's Theorem, see \cite[Theorem 2.1]{nocedalWright99}, it holds that
 $\hat{J}_N(U_0,x_0;\delta)=\hat{J}_N(\hat{U}^*(x_0),x_0;\delta)+\mu_0\nabla_U \hat{J}_N(\hat{U}^*(x_0),x_0;\delta)^\TRANSP d + $ $\frac{1}{2} \mu_0^2 d^\TRANSP \nabla^2_{U} \hat{J}_N(\hat{U}^*(x_0)+sd,x_0;\delta) d$ for some $s\in(0,1)$. 
 But, together with $\hat{J}_N(\hat{U}^*(x_0),x_0;\delta)=\hat{J}_N^*(x_0;\delta)\leq\tilde{J}_N^*(x_0)$, this implies 
 \begin{equation}\label{eq:taylorBound}
  \hat{J}_N(U_0,x_0;\delta) \leq \tilde{J}^*_N(x_0) + \frac{\mu_0^2}{2}\|H+\frac{\varepsilon(1+\bar{w})}{\delta^2} G^\TRANSP G\|_2\, ,
 \end{equation}  
where we used that $\|d\|\leq 1$, $\nabla_U \hat{J}_N(\hat{U}^*(x_0),x_0;\delta)=0$, and $\nabla^2_{U} \hat{J}_N(U,x;\delta) \leq H+\frac{\varepsilon(1+\bar{w})}{\delta^2} G^\TRANSP G$ for any $(U,x) \in \mathbb{R}^{Nm}\times\mathbb{R}^n$, see the proof of Lemma~\ref{lem:jProp}. Let us now choose $\mu_0=\delta$. 
Then, it follows from~(\ref{eq:taylorBound}) that for any $x_0\in\mathcal{X}_0$, any $U_0 \in \mathcal{B}_{\delta}^{Nm}(\hat{U}^*(x_0))$, and any $\delta \leq \delta_0(x_0)$,
\begin{equation} \label{eq:finalBound}
   \hat{J}_N(U_0,x_0;\delta) \leq \tilde{J}^*_N(x_0)+ \frac{1}{2} \|\delta^2 H+\varepsilon(1+\bar{w}) G^\TRANSP G\|\, .
\end{equation} 
 This shows that the cost function value, and hence also $\hat{\alpha}(U_0,x_0;\delta)$ in~(\ref{eq:zNdeltaEquivalence}), will stay bounded for any $\delta \leq \delta_0(x_0)$. \\
 Let us now choose $\delta'_0(x_0)\leq \delta_0(x_0)$ in such a way that the right hand side of~(\ref{eq:zNdeltaEquivalence}) is satisfied for all $\delta \leq \delta'_0(x_0)$, which then implies via (\ref{eq:zNdeltaEquivalence}) that $(U_0,x_0)\in \hat{\mathcal{Z}}_N(\delta)$. 
 In particular, we define $\delta'_0(x_0):=\max\{\delta \in\mathbb{R}_{\pplus}\,|\,\delta\leq\delta_0(x_0),$ $\tilde{J}^*_N(x_0) + \frac{1}{2} \|\delta^2 H+\varepsilon(1+\bar{w}) G^\TRANSP G\| -x_0^\TRANSP P_{\mathrm{uc}}^*\, x_0^{} \leq \varepsilon \bar{\beta}(\delta)\}$.
 In combination with~(\ref{eq:alphaEq}), it follows from~(\ref{eq:finalBound}) that the right hand side of (\ref{eq:zNdeltaEquivalence}) is satisfied for any $\delta\leq\delta'_0(x_0)$. Hence, $(U_0,x_0)\in \hat{\mathcal{Z}}_N(\delta)$ for any $\delta\leq\delta'_0(x_0)$.
 We may picture the construction of $\delta'_0(x_0)$ as decreasing $\delta$ from $\delta(x_0)$ until also the second condition in the foregoing definition is satisfied.
As $\bar{\beta}(\delta)$ is strictly monotonic and grows without bound for decreasing $\delta$, this can always be achieved, which shows that $\delta'(x_0)$ exists for any $x_0\in\mathcal{X}_N^\circ$.
Moreover, since $\bar{\beta}(\delta)$ is continuous in $\delta$ and both $\delta_0(x_0)$ and $\tilde{J}^*_N(x_0)$ are continuous in $x_0$, see Lemmas~\ref{lem:relBarrierLevelSets} and \ref{lem:continuityRes} in the Appendix, ${\delta}'_0(x_0)$ is continuous. 
 \\
\emph{b)\ }
Consider now an arbitrary compact set $\mathcal{X}_0 \subseteq \mathcal{X}_N^\circ$ and define $\bar{\delta}_0\in \mathbb{R}_{\pplus}$ as $\bar{\delta}_0=\min_{x_0\in\mathcal{X}_0}{\delta}'_0(x_0)$. Due to the continuity of
 ${\delta}'_0(x_0)$ and the compactness of $\mathcal{X}_0$, this value always exists by virtue of the Weierstra{\ss} extreme value theorem. 
Then, by the above arguments, $(U_0,x_0)\in \hat{\mathcal{Z}}_N(\delta)$ for any $\delta\leq \bar{\delta}_0$ and any $x_0\in\mathcal{X}_0$ if we choose $U_0$ as $U_0 \in \mathcal{B}_{\mu(\delta)}^{Nm}(\hat{U}^*(x_0))$ with $\mu(\cdot)=\mathrm{id} \in \mathcal{K}_\infty$, which proves our second claim. 
As mentioned above, the zero tolerance case comprises all other cases with nonzero tolerances $\hat{z}_{\mathrm{x},\text{tol}} \in \mathbb{R}^{q_x}_+$, $\hat{z}_{\mathrm{u},\text{tol}} \in \mathbb{R}^{q_u}_+$. \hfill $\blacksquare$
 \begin{rem}
  Note that the above results, including the restriction of $x_0$ to the set $\mathcal{X}_N^\circ$, are more of a conceptual nature and tend to be rather conservative when compared to the actual closed-loop behavior, see also Section~\ref{sec:numExample}. For practical applications, more appropriate choices or selection regimes for $\delta$ and $U_0$ may for example be found based on numerical simulations. 
 \end{rem}
\subsection{Choosing the optimizer update operator} \label{subsec:optOperator}
Apart from condition~(\ref{eq:optOperatorCond}), we provided up to now no explicit characterization of the optimizer update operator $\Psi_{\mathrm{o}}: \mathbb{R}^{Nm}\times \mathbb{R}^n \to \mathbb{R}^{Nm}$. In the following, we propose to choose $\Psi_{\mathrm{o}}(\cdot,\cdot)$ based on a line search method of the form
\begin{equation} \label{eq:optOperatorP}
\Psi_{\mathrm{o}}(U,x)=U+s\, p(U,x)\, ,
\end{equation}
where $s \in\mathbb{R}_{\pplus}$ is a step size parameter and $p:\mathbb{R}^{Nm}\times \mathbb{R}^n \to \mathbb{R}^{Nm}$ is chosen according to a a suitable search direction rule. Note that the step size~$s$ is not fixed but has to be chosen at each iteration step based on the current search direction.
In order to ensure a cost function decrease at each optimizer update, we are especially interested in search direction and step size selection approaches which ensure that $p$ and $s$ satisfy the so-called 
strong Wolfe conditions~\cite{nocedalWright99}
\begin{subequations} \label{eq:wolfeCond}
 \begin{align}
   \hat{J}_N(U+sp,x)\leq \hat{J}_N(U,x)+c_1 s\, \nabla_{\!U} \hat{J}_N(U,x)^{\!\TRANSP}\! p \\
   |\nabla_U \hat{J}_N(U+sp,x)^\TRANSP p|\leq c_2|\nabla_U \hat{J}_N(U,x)^\TRANSP p|
 \end{align}
\end{subequations} 
for positive constants $c_1, c_2 \in \mathbb{R}_{\pplus}$ with $0<c_1<c_2<1$.
In fact, (\ref{eq:wolfeCond}a) ensures that~(\ref{eq:optOperatorCond}) will hold with $\gamma(U,x)=-c_1s\,\nabla_{\!U} \hat{J}_N(U,x)^{\!\TRANSP} p(U,x)$ if $p=p(U,x)$ is chosen as a \emph{descent direction}, i.e., if $\nabla_{\!U} \hat{J}_N(U,x)^{\!\TRANSP}\! p(U,x) <0$ whenever $\nabla_{\!U} \hat{J}_N(U,x) \neq 0$. 
Suitable approaches for the choice of $p=p(U,x)$ are for example given by
the gradient method (G), the conjugated gradient method (CG), the Newton method~(N), or the Quasi-Newton method~(QN), i.e.
\begin{align}
G\!: \ p(U,x)& =-\nabla_U \hat{J}_N(U,x) \, ,  \label{eq:searchGradient} \\[0.15cm]
CG\!:\ p(U,x)& =-\nabla_U \hat{J}_N(U,x) +\beta_{\mathrm{CG}}\, p_{-1}\, , \label{eq:searchConjGradient} \\[0.15cm]
N\!: \ p(U,x)& =-\left(\nabla^2_U \hat{J}_N(U,x)\right)^{\!-1}\!\nabla_U \hat{J}_N(U,x) \, , \label{eq:searchNewton} \\[0.15cm]
QN\!:\ p(U,x)& =-\hat{B}\, \nabla_U \hat{J}_N(U,x)\, \label{eq:searchQNewton}\, . 
\end{align}
Here, $\hat{B} \in\mathbb{S}^{Nm}_{\pplus}$ is a suitably updated approximation of the inverse Hessian matrix, for example computed by using the well-known BFGS formula, while $\beta_{\mathrm{CG}}$ in~(\ref{eq:searchConjGradient}) is a scalar that ensures that $p(U,x)$ is conjugate to the previously used search direction $p_{-1}$, see~\cite{nocedalWright99} for more details. 
In combination with a suitable step size selection procedure, all the aforementioned approaches can be shown to provide a descent direction, which implies that any of them may be used within the optimizer update operator of the above MPC iteration scheme. 
In particular, we can state the following corollary.
%
\begin{cor}\label{cor:optOperatorRes}
Let the cost function $\hat{J}_N:\mathbb{R}^{Nm}\times\mathbb{R}^n \to \mathbb{R}_+$ given in~(\ref{eq:QPformulation}) be associated to the relaxed barrier function based MPC approach based on~(\ref{eq:relBarrierDef})--(\ref{eq:riccatiP}).  
Moreover, consider the closed-loop system~(\ref{eq:algoDyn}) with the operator $\Phi^{i_{\mathrm{T}}(k)}$ defined according to~(\ref{eq:recIteration}), (\ref{eq:shiftOperator}), and (\ref{eq:optOperatorP}), where $p=p(U,x)$ is chosen according to one of the approaches in (\ref{eq:searchGradient})--(\ref{eq:searchQNewton}) with $p$ and $s$ satisfying the strong Wolfe conditions~(\ref{eq:wolfeCond}) at each internal iteration. Then, for any sequence $\boldsymbol{i}_{\mathrm{T}}=\{i_{\mathrm{T}}(0), i_{\mathrm{T}}(1), \dots\}$, the origin $(U,x)=(0,0)$ of system~(\ref{eq:algoDyn}) is globally asymptotically stable. 
Moreover, concerning the maximal possible constraint violations, Theorem~\ref{thm:maxConstrViolation} and~ Theorem~\ref{thm:xNdelta0} apply. 
\end{cor}
As the cost function~$\hat{J}_N(\cdot,x)$ is continuously differentiable and bounded from below and $p=p(U,x)$ is a descent direction, there always exists a nonempty interval of step lengths such that the strong Wolfe condition~(\ref{eq:wolfeCond}) will be satisfied, see Lemma~3.1 in~\cite{nocedalWright99}. 
For a discussion of suitable step size selection procedures we refer the interested reader to~\cite[Chapter 3.4]{nocedalWright99}.\\
A particular interesting case is to use the Newton search direction from~(\ref{eq:searchNewton}) together with a so-called backtracking line search, in which a suitable step size is found by increasing an integer parameter $j$ starting from $j=0$ until $s=\rho^{j}$ satisfies the Armijo condition~(\ref{eq:wolfeCond}a), where $\rho\in (0,1)$ is an a priori chosen design parameter. With the help of~(\ref{eq:Jprop}a), it can be shown that~(\ref{eq:wolfeCond}a) will be satisfied for all $s\leq \bar{s}:=2 \sigma (1-c_1)/L$, see Appendix~\ref{sec:appStepSize}, which implies that the backtracking line search will always terminate with a step size not smaller than $2 \rho  \sigma (1-c_1)/L$. From this, we can derive the following upper bound on the number of backtracking line search iterations
  \begin{equation}
 j \leq {j}_{\mathrm{max}}:= 1+ \log_{\rho}\left(2\sigma(1-c_1)/L\right) \, ,
 \end{equation}
which then yields the following upper bound for the overall complexity that is inferred by the proposed MPC iteration scheme in each sampling step
\begin{equation}
C_{\mathrm{iter}}(k)\leq C_{\mathrm{shift}}+i_{\mathrm{T}}(k)\left({j}_{\mathrm{max}}\, C_{\mathrm{bt}} + C_{\mathrm{nd}} \right) \, .
 \end{equation} 
Here, $C_{\mathrm{shift}}$ and $C_{\mathrm{bt}}$ refer to the computational complexities of performing one shift update and one backtracking line search iteration, respectively, while~$C_{\mathrm{nd}}$ denotes the complexity of computing the Newton direction. It can be expected that the complexity of the Newton step will in most cases dominate the other terms. A thorough complexity analysis as well as a discussion on how the complexity of the optimizer update~(\ref{eq:optOperatorP}) can be reduced is the topic of ongoing research.
Nevertheless, the above arguments already allow to compute a rough estimate of the number of possible algorithm iterations for a given sampling time, or, vice versa, to choose a suitable sampling time for a desired number of optimizer updates.
\section{Numerical Example}\label{sec:numExample}
\begin{figure*}[ht]
\hspace*{-0.15cm}\includegraphics[scale=0.45]{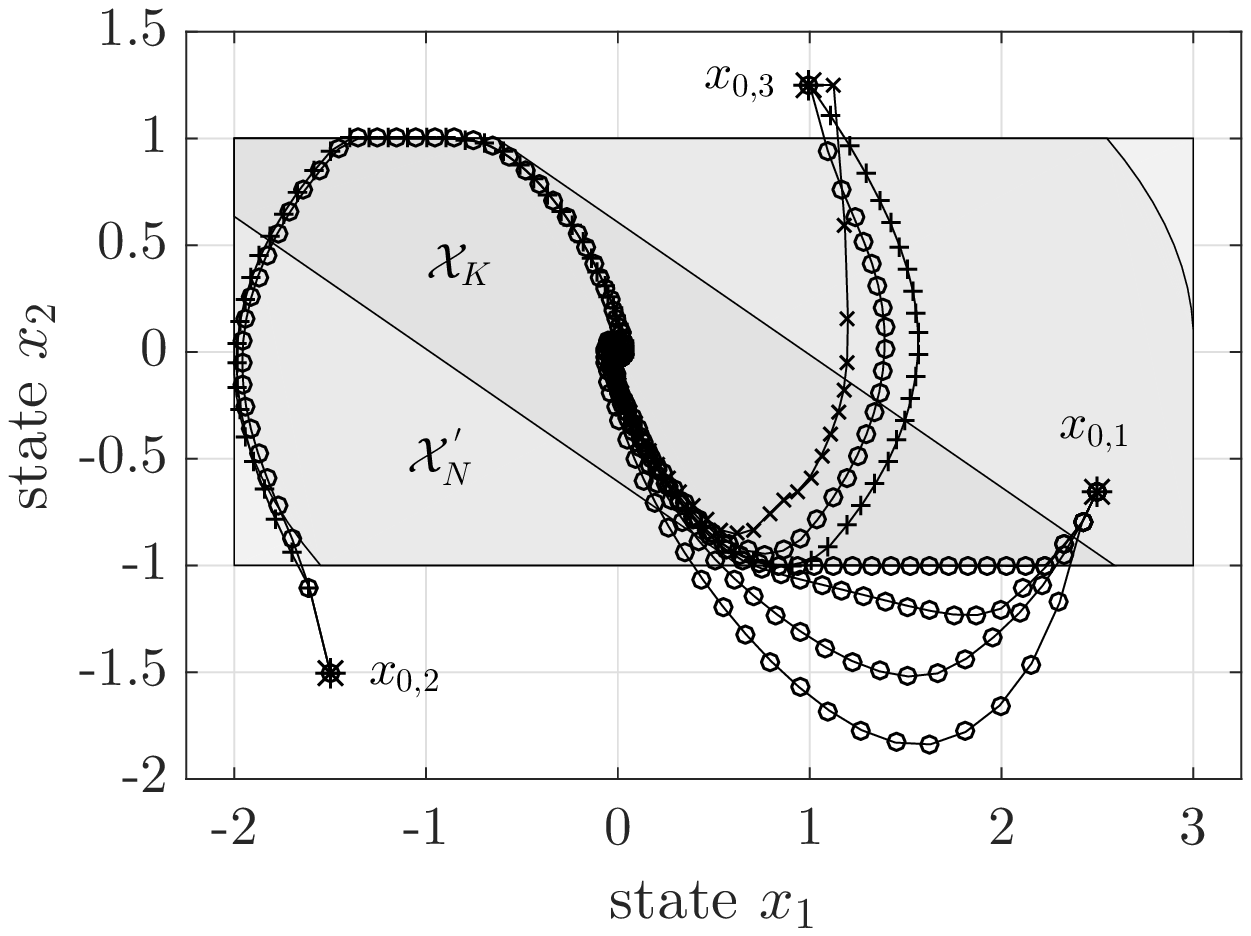}
\hspace*{-0.25cm}\includegraphics[scale=0.45]{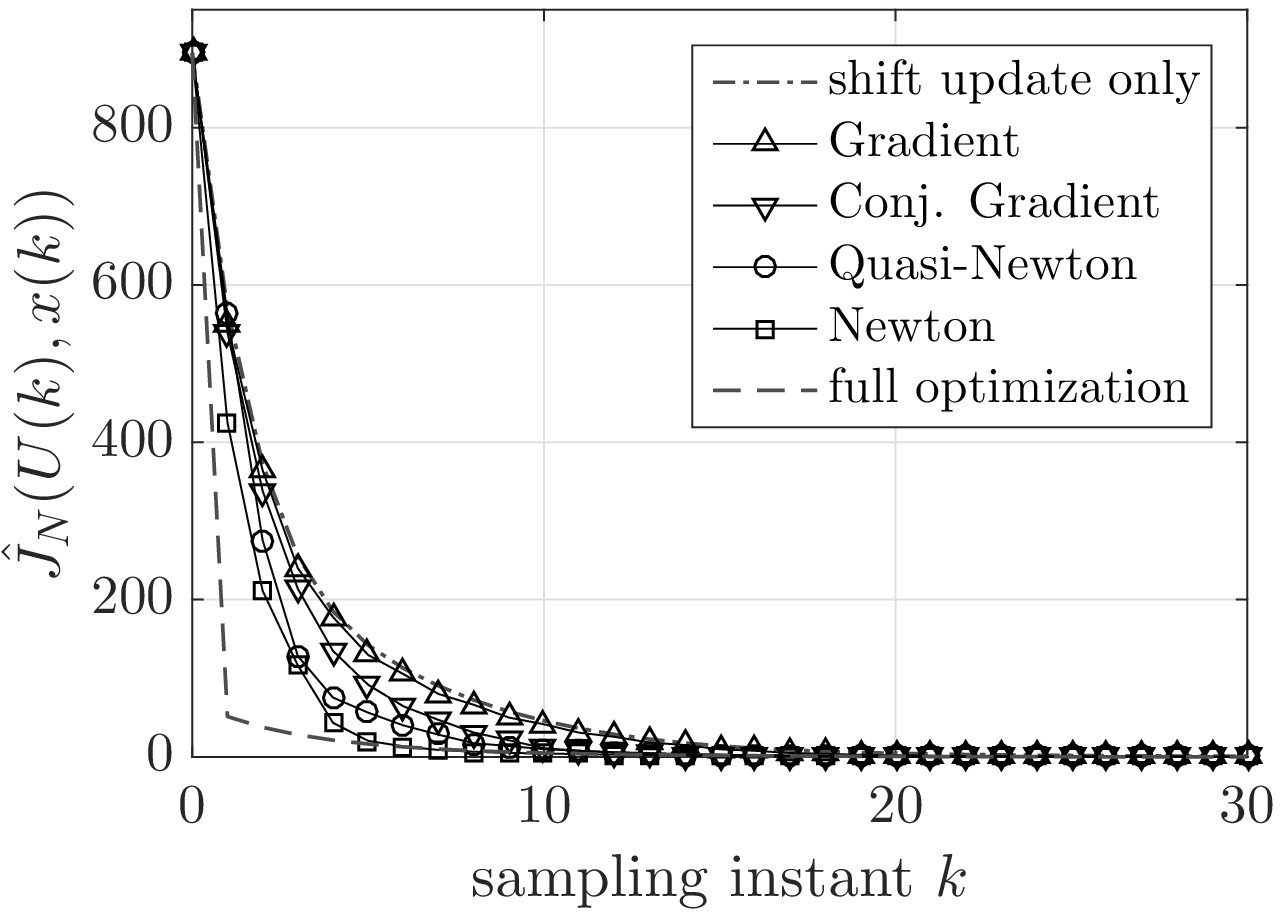}
\hspace*{-0.25cm}\includegraphics[scale=0.45]{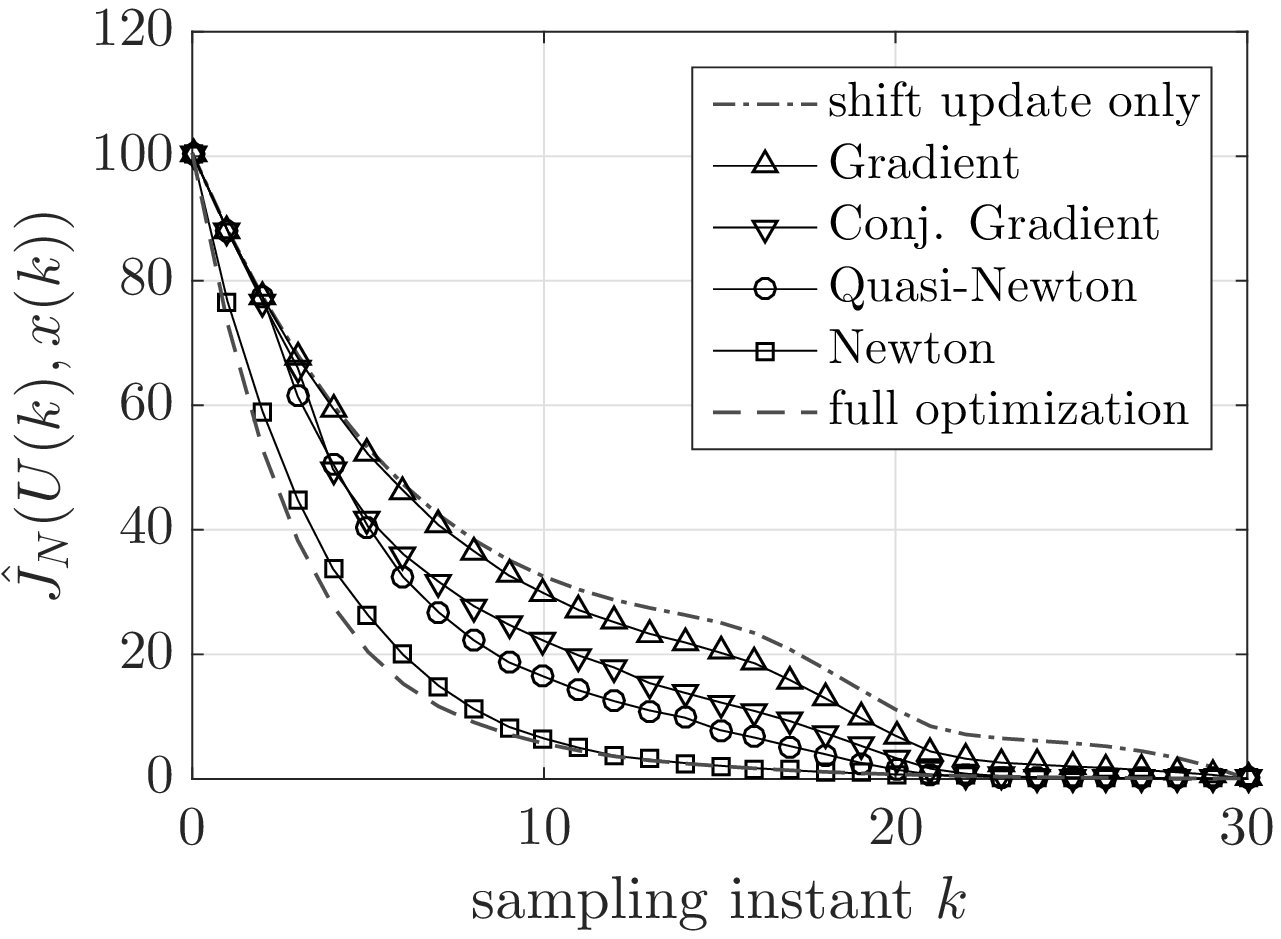}
\caption{\emph{Left:} Closed-loop state trajectories for the rbMPC iteration scheme when applying only one Newton iteration per sampling step. For $x_{0,1}=[2.5,-0.65]$, we vary the relaxation parameter $\delta \in \{0.1,2.5\times10^{-2},10^{-2},10^{-3}\}$; for $x_{0,2}=[-1.5,-1.5]$, we plot for comparison also the optimal trajectory $[+]$ corresponding to $i_{\mathrm{T}}(k)=100$; for $x_{0,3}=[1,1.25]$, we compare the behavior for the initializations $U_0=\bar{K}x_0$~$[\circ]$, $U_0=\hat{U}^*(x_0)$~$[+]$, and $U_0=0$~$[\times]$. Also depicted are the constraint set $\mathcal{X}$, the set $\mathcal{X}_K=\{x\in\mathcal{X}:Kx\in\mathcal{U}\}$ and the feasible set $\mathcal{X}_N'=\{x\in\mathcal{X}: \exists U\in \mathbb{R}^{Nm} \text{ s.t. } GU\leq Ex+d\}$.  \newline 
\emph{Middle:} Averaged evolution of the cost function for 200 random initial conditions $x_0$ and $U_0=\bar{K}x_0$. Newton based optimization achieves the fastest decay, while the gradient method is almost similar to applying only the shift operator ($i_{\mathrm{T}}(k)\equiv 0$). \newline  
 \emph{Right:} Averaged evolution of the cost function for 200 random initial conditions $x_0$ and $U_0=\hat{U}^*(x_0)$. Performing a single Newton step at each sampling instant achieves almost the same behavior as when applying the optimal solution.
}
\label{fig:numExample}
\end{figure*} 
We want to briefly illustrate some aspects of the discussed barrier function based MPC iteration scheme by means of an academic numerical example.
The discussion deliberately focuses on the systems theoretic properties analyzed above, leaving a thorough investigation of the related numerical aspects for future work.   
Considered is a double integrator system of the form
 \begin{equation}
x(k+1)=\begin{bmatrix}
1 & T_s\\
0 & 1
\end{bmatrix}x(k) + \begin{bmatrix}
T_s^2 \\ T_s
\end{bmatrix} u(k)\,, \ \ T_s=0.1\, ,
\end{equation}
with the input and state constraint sets $\mathcal{U}$ and $\mathcal{X}$ being defined by $\mathcal{U}=\{u \in \mathbb{R}: |u| \leq 1\}$ and $\mathcal{X}=\{x \in \mathbb{R}^2: -2\leq x_1 \leq 3, |x_2| \leq 1\}$, respectively.
Prediction horizon, weight matrices, and barrier function parameters are chosen as $N=30$, $Q=\mathrm{diag}(1, \ \, 0.1)$, $R=0.1$, $\varepsilon=10^{-3}$, and $\delta=10^{-3}$ if not stated otherwise. The weighting vectors $w_{\mathrm{x}}\in\mathbb{R}^{q_x}_+$, $w_{\mathrm{u}}\in\mathbb{R}^{q_u}_+$ for the recentering are computed according to the procedure presented in~\cite{feller15b}. The barrier function based MPC cost function is formulated based on the design from Section~\ref{sec:relBarrierIntro}, and the proposed iteration scheme is, combined with the line search approaches from Section~\ref{subsec:optOperator}, implemented and tested for various initializations $(U_0,x_0)$. A backtracking line search with $\rho=0.5$ is used for the gradient and Newton method, while the more involved step size procedure from Mor\'{e} and Thuente is used for both the conjugated gradient and Quasi-Newton method, see~\cite{more94}. The Wolfe condition parameters are chosen as $c_1=10^{-3}$, $c_2=0.9$. The inverse Hessian approximation $\hat{B}$ in the Quasi-Newton method is initialized with $H^{-1}$ and then handed over between consecutive sampling instances. Fig.~\ref{fig:numExample} summarizes exemplary closed-loop results for the case of applying only \emph{one} optimizer update per sampling step, i.e. $i_{\mathrm{T}}(k)=1 \ \forall\, k \in \mathbb{N}$. \\
The left subplot of Fig.~\ref{fig:numExample} illustrates that a quite good closed-loop performance, including the satisfaction of state constraints, can be achieved even when the iteration scheme is initialized with a suboptimal (and in many cases infeasible) $U_0=\bar{K} x_0$. For the example at hand,
\begin{equation}
 \bar{K}=[K^\TRANSP\ \ \dots \ {(A+BK)^{N-1}}^\TRANSP K^\TRANSP] \, ,
\end{equation} 
which corresponds to initialization with the open-loop input sequence resulting from the controller gain $K$ in~(\ref{eq:riccatiP}a). It turns out that, although the initialization is far from optimal, almost exact constraint satisfaction (in input and states) can be achieved for $x_{0,1}$ whenever $\delta\leq 10^{-3}$. This illustrates the practical relevance, but, to a certain extent, also the conservatism of the results presented in Section~\ref{subsec:constraints}. Moreover, it can be seen that the approach even allows for infeasible initial conditions $x_0 \notin \mathcal{X}$ and that, at least in the case of a Newton based optimization, neither the number of performed optimizer updates nor the quality of the initialization $U_0$ has a very big impact on the resulting closed-loop behavior.  \\
The middle and right subplots of Fig.~\ref{fig:numExample} show, averaged over 200 randomly chosen initial conditions $x_0$ from the set $\mathcal{X}$, the strictly monotonic decay of the cost function over time when applying the different line search methods discussed above, based on the initializations $U_0=\bar{K}x_0$ and $U_0=\hat{U}^*(x_0)$, respectively. Again, only one optimizer update is performed per sampling step. It can be seen that particularly the Newton method based optimization does result in a good closed-loop performance, especially in the case of optimal initialization. 
The observed behavior is in principle also reflected on the level of maximal constraint violations. In particular, while one Newton iteration is often enough to achieve almost the same constraint satisfaction properties as the fully optimized solution, the Quasi-Newton and conjugated gradient based methods tend to require a few more, say five to ten, iterations. 
This reflects the fact that these approaches may need some internal iterations in order to reconstruct the second order information that is in the Newton method directly provided by the Hessian. 
The gradient method performs rather poorly, which, in principle, had to be expected as the considered optimization problems are often ill-conditioned due to the underlying relaxation. \\
Thus, although the resulting closed-loop behavior may depend on the line search method that is used within the optimizer update, we are, in summary, able to reproduce and confirm many of the systems theoretic properties discussed above. In particular, asymptotic stability of the overall closed-loop system is observed for all initial conditions and independently of both the applied optimization procedure and the number of optimizer iterations.      
Our simulations furthermore suggest that the Newton based line search often converges already after only five to eight iterations, which nicely relates to the observations reported in~\cite{wang10}. 
On the other hand, the above results also show a promising behavior of the optimizer update approaches based on a Quasi-Newton or conjugated gradient line search, which suggests that knowledge and inversion of the full Hessian matrix might not always be necessary.  
\section{Conclusion}\label{sec:conclusion}
We presented and analyzed an MPC iteration scheme which is based on the concept of relaxed logarithmic barrier functions and performs only a limited number of optimization algorithm iteration between two consecutive sampling instants. 
The stability properties of the overall closed-loop system, consisting of system state and optimization algorithm dynamics, were discussed and we showed that, under standard assumptions on the optimization procedure and the system to be controlled, asymptotic stability of the origin can be guaranteed independently of the number of optimization algorithm iterations.
In addition, we studied the corresponding constraint satisfaction properties and showed that, assuming an appropriate initialization of the optimization algorithm, both approximate and exact constraint satisfaction may in many cases be recovered by a suitable design of the underlying relaxed barrier functions.  
\newpage
\noindent These results were also illustrated by means of the presented numerical example. \\
As outlined above, the discussed iteration scheme can be seen as a stabilizing \emph{anytime MPC algorithm} in the sense that it always returns an asymptotically stabilizing control input and improves the quality of the returned control at each of its internal iterations. 
This shows the potential relevance of the presented approach for practical applications in which the time and hardware resources that are available for the on-line optimization may be limited or even varying over time.  \\
Interesting open problems include a robustness analysis of the discussed iteration scheme, a thorough investigation of the numerical implementation and performance, as well as the design of tailored optimization algorithms that may use properties or structure of the underlying barrier function based formulation. 
Finally, some of the above results may not necessarily be limited to the discussed barrier function based approach but could also be applicable to existing MPC schemes based on penalty functions or soft constraints, as well as to more general classes of parametric optimization problems not necessarily related to~MPC.
\appendix
 \section*{Appendix}
 \renewcommand{\thesubsection}{A.\arabic{subsection}}
%
\subsection{Condensed problem formulation} \label{app:condensedForm}
Note that for given $x\in\mathbb{R}^n$ and $U\in\mathbb{R}^{Nm}$, the predicted system states for $k=1,\dots,N$ can be expressed as
\begin{equation}\label{eq:xkDyn}
 x_k(U,x)=A^kx + \sum_{i=0}^{k-1} A^iBu_{k-1-i} \, .
\end{equation}
By writing~(\ref{eq:OptProblemRelBarrier}) in matrix form and eliminating the predicted system states by means of~(\ref{eq:xkDyn}), we get that the cost function matrices in~(\ref{eq:QPformulation}) are given by
\begin{subequations}
 \begin{align*}
  H&=2\left(\tilde{R}+ \Gamma^{\!\TRANSP}\! \tilde{Q} \Gamma \right),  F=2\Omega^{\!\TRANSP}\! \tilde{Q}\Gamma, Y=2\left({Q}+ \Omega^{\!\TRANSP}\! \tilde{Q} \Omega \right) \\
  \Omega&=\begin{bmatrix}A \\[-0.1cm] \vdots \\ A^N \end{bmatrix},
  \Gamma=\begin{bmatrix} B & \cdots & 0 \\[-0.1cm]
  \vdots & \ddots & \vdots \\
  A^{N-1}B  & \cdots & B
  \end{bmatrix}, \begin{array}{l} 
                  \tilde{Q}= I_N \otimes Q \\[0.3cm]
                   \tilde{R}=I_N \otimes R\, ,
                \end{array}
 \end{align*}
 \end{subequations}
 while the constraint matrices $G$, $d$, and $E$ can be constructed as
 \begin{subequations}
  \begin{align*}
    G&=\begin{bmatrix} 0 \\ ([I_{N-1}\ \,0] \otimes C_{\mathrm{x}})\Gamma \\ I_N \otimes C_{\mathrm{u}}\end{bmatrix}, \
 E=\begin{bmatrix} I_n \\([I_{N-1}\ \,0] \otimes C_{\mathrm{x}}) \Omega \\ 0 \end{bmatrix} \\
  d&=\begin{bmatrix} d_{\mathrm{x}}^\TRANSP & \cdots & d_{\mathrm{u}}^\TRANSP \end{bmatrix}^\TRANSP,
  \end{align*}
 \end{subequations}
Note that the representation of $G$, $d$, and $E$ may change if redundant constraints are removed by making use of a suitable algorithm.
\subsection{Proof of Lemma \ref{lem:jProp}} \label{app:barrierProperties}
Due to compactness of the constraint sets $\mathcal{X}$ and $\mathcal{U}$ and the design of the corresponding relaxed barrier functions, both $\hat{\ell}(x,u)$ and $\hat{F}(x)$ are continuous, positive definite, and radially unbounded. Since $x_k$, $k=1,\dots,N$, depends linearly on $x$ and $U$ via~(\ref{eq:xkDyn}), the same holds for $\hat{J}_N(U,x)$, which implies that $\underline{\alpha}$, $\bar{\alpha} \in \mathcal{K}_\infty$ satisfying property~$i)$ will always exist, see Lemma~4.3 in \cite{khalil02}.\\
Strong convexity in $U$ and convexity in~$x$ follow from $R\in \mathbb{S}^m_{\pplus}$, $Q\in\mathbb{S}^n_+$ and convexity of the relaxed barrier functions. In fact, the overall barrier function reads
\begin{equation}
\hat{B}_{\mathrm{xu}}(U,x)=\!\sum_{i=1}^{q} \left(1+w^i\right)\!\big(\hat{B}\left(z_i(U,x)\right) - \ln\left(d^i\right)\big),
\end{equation} 
where the relaxed logarithmic barrier function~$\hat{B}(\cdot)$ is defined according to (\ref{eq:relBarrierDef}) with $z_i(U,x)=-G^iU+E^ix+d^i$, while the weighting vector $w \in \mathbb{R}^q_+$ is used for recentering around the origin.
Based on this, we can easily compute the Hessian of the cost function with respect to $U$ as 
\begin{equation}\label{eq:hessian}
 \nabla^2_U \hat{J}_N(U,x)=H + \varepsilon G^\TRANSP \mathrm{diag}\big(D_{1}(U,x), \dots, D_{q}(U,x)\big) G \, ,
\end{equation}
with $D_i(U,x)$ for $i=1,\dots,q$ being defined as
\begin{equation}
 D_i(U,x)=\begin{cases} \frac{1+w^i}{z_i(U,x)^2} & \, z_i(U,x) > \delta \\
 \ \ \frac{1+w^i}{\delta^2} & \, z_i(U,x) \leq \delta \, . \end{cases}
\end{equation}
Obviously, it holds that $0\leq D_i(U,x) \leq \frac{1+\bar{w}}{\delta^2}$ with $\bar{w}:=\max_i \{w^i\}$, which immediately yields that 
\begin{equation}\label{eq:upperLowerBounds}
 \sigma I \preceq \nabla^2_U \hat{J}(U,x) \preceq LI \, ,
\end{equation} 
with $\sigma=\lambda_{\min}(H)$, $L=\lambda_{\max}(H+\frac{\varepsilon(1+\bar{w})}{\delta^2}G^\TRANSP G) \in \mathbb{R}_{\pplus}$. Thus, the cost function is not only strongly convex but we can also give a global upper bound on the Hessian, and (\ref{eq:Jprop}a) holds. 
By means of the mean value theorem for multivariable functions we further obtain
\begin{align}
 \nabla_U \hat{J}_N&(U_1,x)-\nabla_U \hat{J}_N(U_2,x) = \\ 
 &\int_{0}^1 \nabla_U^2 \hat{J}_N(U_2+s(U_1-U_2),x)\,\mathrm{d}s \ (U_1-U_2) \nonumber
 \end{align}
and hence, with the above bounds on the Hessian,
\begin{equation} \label{eq:jLipschitzBoundProof}
\|\nabla_U\hat{J}_N(U_1,x)-\nabla_U\hat{J}_N(U_2,x)\| \leq L\|U_1-U_2\| \, 
\end{equation} 
for all $U_1, U_2 \in \mathbb{R}^{Nm}$, $x\in \mathbb{R}^n$, which proves property~$ii)$. \\
A proof for property $iii)$ has been presented in the proof of Theorem~5 in~\cite{feller15a} and only a sketch is given here. Based on the quadratic relaxation, the barrier functions can be shown to satisfy $\hat{B}_{\mathrm{x}}(x)\leq x^{\TRANSP} M_{\mathrm{x}}\, x^{} \ \forall\, x\in \mathbb{R}^n$ as well as $\hat{B}_{\mathrm{u}}(u)\leq u^{\TRANSP} M_{\mathrm{u}}\, u^{} \ \forall\, u\in \mathbb{R}^m$, see Lemma~3 in \cite{feller15a}. Let now $K$ and $P$ be given by~(\ref{eq:riccatiP}) and choose $k_f(U,x)=K x_N(U,x)$, where $x_N(U,x)$ is given according to~(\ref{eq:xkDyn}). Then, the aforementioned quadratic upper bounds ensure that $\hat{F}((A+BK)x_N)-\hat{F}(x_N) \leq -\hat{\ell}(x_N,Kx_N) \ \forall\, x_N \in \mathbb{R}^n$. Defining $U_0^+$ as in Lemma~\ref{lem:jProp} and using a telescoping sum, this directly implies that~(\ref{eq:jPropDecrease}) holds. Thus, property~$iii)$ holds with $k_f(U,x)=K x_N(U,x)$. $\hfill \blacksquare$
\subsection{Some auxiliary continuity and feasibility results} \label{app:auxResults}
\begin{lem} \label{lem:relBarrierLevelSets}
Let $\hat{B}_{\mathcal{P}}(\cdot)$ be a recentered relaxed logarithmic barrier functions of the form~(\ref{eq:relBarrierPol}) for a polytopic set $\mathcal{P}=\{x\in\mathbb{R}^r : C \xi \leq d \}$ with $C\in\mathbb{R}^{q\times r}$ and $d\in\mathbb{R}^q_{\pplus}$. Let $\bar{\beta}(\delta)$ be defined as $\bar{\beta}(\delta):=\min_{i,\xi}\{\hat{B}_{\mathcal{P}}(\xi)\vert \ {C^i\xi=d^i}\}$ with $i=1,\dots,q$. Then, $\bar{\beta}(\cdot)$ is a continuous function. Furthermore, $\mathcal{S}_{\hat{B}}=\{x \in \mathbb{R}^n| \hat{B}_{\mathcal{P}}(x)\leq \bar{\beta}(\delta)\}\subseteq \mathcal{P}$ for any $\delta \in \mathbb{R}_{\pplus}$. 
\end{lem}
\textit{\ \ Proof: }
The barrier function~$\hat{B}_{\mathcal{P}}(\cdot)$ is continuous in both $\xi$ and $\delta$ and possesses compact sublevel sets for any fixed value of~$\delta$. This implies that $\bar{\beta}_i(\delta):=\min_{\xi}\{\hat{B}_{\mathcal{P}}(\xi)\vert \ {C^i\xi=d^i}\}$ will be continuous in $\delta$ for any $i=1,\dots,q$, see~Proposition~4.4 in~\cite{bonnans00}. Hence, $\bar{\beta}(\delta):=\min_{i} \bar{\beta}_i(\delta)$ will also be continuous. \\
Assume there exists an $\bar{\xi}\in \mathcal{S}_{\hat{B}}$ with $\bar{\xi} \notin \mathcal{P}$, i.e., $\hat{B}_{\mathcal{P}}(\bar{\xi})\leq\bar{\beta}(\delta)$ and $C \bar{\xi}>d$. However, then there would exists a $\lambda \in (0,1)$ such that if $\xi:=\lambda \bar{\xi}$, then $C \xi \leq d$ and $C^i \xi=d^i$ for some $i=1,\dots,q$. Now, due to the convexity and positive definiteness of $\hat{B}_{\mathcal{P}}(\cdot)$, it holds that $\hat{B}_{\mathcal{P}}(\xi)=\hat{B}_{\mathcal{P}}(\lambda\bar{\xi})\leq \lambda \hat{B}_{\mathcal{P}}(\bar{\xi})<\bar{\beta}(\delta)$, which is a contradiction to the definition of $\bar{\beta}(\delta)$ as $C^i \xi=d^i$. Thus, $\bar{\xi}\in \mathcal{P}$ whenever $\bar{\xi}\in\mathcal{S}_{\hat{B}}$, which completes the proof. $\hfill \blacksquare$ \\[-0.2cm]
\begin{lem} \label{lem:continuityRes}
  Assume $(A,B)$ to be controllable and let $\mathcal{X}_N:=\{x\!\in\!\mathbb{R}^n\!:\!\exists\, U\!\in \mathbb{R}^{Nm} \mbox{\ s.t.\ } u_k \in \mathcal{U},$ $x_k(U,x)\in \mathcal{X} \mbox{\ for \ } k=1,\dots,N-1, x_N(U,x)=0 \}$ be the set of states for which the problem
 \begin{subequations}\label{eq:OptProblemBarrier}
\begin{align}
\tilde{J}_N^{\ast}(x)&=\min_{U}\sum_{k=0}^{N-1} \tilde{\ell}(x_k,u_k)   \\
\mbox{s.\,t.} \ \ & {x}_{k+1}=A{x}_{k}+B{u}_{k}, \ x_N=0, \ {x}_0=x\, .
\end{align}
\end{subequations}
has a solution. Here, $\tilde{\ell}(x,u):= \|x\|_Q^2+\|u\|_R^2 + \varepsilon B_{\mathrm{u}}(u)+\varepsilon B_{\mathrm{x}}(x)$, where $B_{\mathrm{u}}(u)$ and $B_{\mathrm{x}}$ are recentered logarithmic barrier functions, e.g. barrier functions of the form~(\ref{eq:relBarrierPol}) with $\hat{B}(\cdot)$ being replaced by $B(z)\!=\!-\ln(z)$. Then, for any $x_0 \in \mathcal{X}_N^\circ$,
the optimizer $\tilde{U}^*(x_0)$ associated to~(\ref{eq:OptProblemBarrier})
is unique and $\tilde{U}^*:\mathcal{X}_N^\circ \to \mathbb{R}^{Nm}$  is a continuously differentiable function.
\end{lem}
\textit{\ \ Proof:}
Using the input vector notation, problem~(\ref{eq:OptProblemBarrier}) is equivalent to 
  \begin{equation}\label{eq:OptProblemBarrier2}
\tilde{J}_N^{\ast}(x)=\min_{U} \tilde{J}_N(U,x) \mbox{\ \ s.\,t.} \ \  S_1 U +S_2x=0
\end{equation}
where $\tilde{J}_N(\cdot,\cdot)$ is the nonrelaxed barrier function based cost function similar to~(\ref{eq:QPformulation}), strongly convex in $U$, and $S_1=\begin{bmatrix} A^{N-1}B & \cdots & B \end{bmatrix}                                                                                                                                                   $, $S_2=A^N$ are used to encode the terminal equality constraint $x_N(U,x)=0$, see~(\ref{eq:xkDyn}). For given $x_0 \in \mathcal{X}_N^\circ$, the KKT conditions for~(\ref{eq:OptProblemBarrier}), which are in this case necessary and sufficient conditions for optimality, read
\begin{subequations}\label{eq:KKTcond}
\begin{align}
  \nabla_U \tilde{J}_N(U^*(x_0),x_0) + S_1^\TRANSP \nu^*(x_0) &=0 \\
  S_1 \tilde{U}^*(x_0) +S_2x_0 &=0 \, ,
\end{align}
\end{subequations}
with $\nu^*(x_0)$ denoting the associated Lagrange multiplier. Now, the implicit function theorem states that unique and continuously differentiable solutions $\tilde{U}^*(x)$, $\nu^*(x)$ to~(\ref{eq:KKTcond}) will exist if the corresponding Jacobian with respect to $(\tilde{U}^*,\nu^*)$, given by 
\begin{equation}
 \begin{bmatrix}
  \nabla_U^2 \tilde{J}_N(U^*(x_0),x_0) & S_1^\TRANSP \\
  S_1 & 0
 \end{bmatrix}\, ,
\end{equation} 
is nonsingular for all feasible $x_0$. Indeed, this follows from applying block inversion and noting that both the upper left corner block and the Schur complement are nonsingular. In particular, $\nabla_U^2 \tilde{J}_N(U^*(x_0),x_0) \in \mathbb{S}^{Nm}_{\pplus}$ $ \forall\, x_0 \in \mathcal{X}_N^\circ$ (strong convexity) and $S_1 \nabla_U^2 \tilde{J}_N(U^*(x_0),x_0) S_1^\TRANSP \in \mathbb{S}^{n}_{\pplus} \ \forall\, x_0 \in \mathcal{X}_N^\circ$ ($S_1$ has full row rank due to controllability of $(A,B)$).
Note that, using this result, it can be shown easily that the value function $\tilde{J}_N^*:\mathbb{R}^n \to \mathbb{R}_+$ is in fact twice continuously differentiable. $\hfill \blacksquare$
%
\subsection{Step size bound for the Newton method} \label{sec:appStepSize}
For the Newton method search direction from~(\ref{eq:searchNewton}), the Armijo condition~(\ref{eq:wolfeCond}a) can be rewritten as 
\begin{equation} \label{eq:armijoNewton}
 \hat{J}_N(U+sp,x)\leq \hat{J}_N(U,x) - c_1 s\, \lambda(U,x)^2\, 
\end{equation} 
with the corresponding squared Newton decrement
\begin{equation}
 \lambda(U,x)^2:=\nabla_U\hat{J}_N(U,x)^{\TRANSP} \left(\nabla^2_U \hat{J}_N(U,x)\right)^{\!-1}\!\nabla_U \hat{J}_N(U,x). \nonumber
\end{equation} 
Applying Taylor's Theorem to the left-hand side of~(\ref{eq:armijoNewton}) and using the upper and lower bounds from~(\ref{eq:Jprop}a), we obtain that 
\begin{equation}
 \hat{J}_N(U+sp,x)\leq \hat{J}_N(U,x)-s\Big(1-\frac{L}{2\sigma}s\Big)\lambda(U,x)^2\!,
\end{equation}
which directly implies that~(\ref{eq:armijoNewton}) will be satisfied for any $s\leq \bar{s}:=2\sigma(1-c_1)/L$.

\bibliographystyle{plain}
\bibliography{ctrbMPC}

\begin{thebibliography}{10}

\bibitem{bemporad15}
A.~Bemporad, D.~Bernardini, and P.~Patrinos.
\newblock A convex feasibility approach to anytime model predictive control.
\newblock {\em arXiv:1502.07974}, 2015.

\bibitem{bemporadExplicitLQR}
A.~Bemporad, M.~Morari, V.~Dua, and E.~N. Pistikopoulos.
\newblock {The explicit linear quadratic regulator for constrained systems}.
\newblock {\em Automatica}, 38:3--20, 2002.

\bibitem{benTal92}
A.~Ben~Tal, M.~Tsibulevskii, and I.~Yusefovich.
\newblock Modified barrier methods for constrained and minimax problems.
\newblock Technical report, Optimization Laboratory, Technion, 1992.

\bibitem{bonnans00}
J.~F. Bonnans and A.~Shapiro.
\newblock {\em Perturbation analysis of optimization problems}.
\newblock Springer, 2000.

\bibitem{diehl07}
M.~Diehl, R.~Findeisen, and F.~Allg\"ower.
\newblock {A stabilizing iteration scheme for nonlinear model predictive
  control}.
\newblock {\em In: Real-Time and Online {PDE}-Constrained Optimization}, pages
  23--52, 2007.

\bibitem{diehl05}
M.~Diehl, R.~Findeisen, F.~Allg\"ower, H.~Bock, and J.~Schl\"oder.
\newblock {Nominal stability of real-time iteration scheme for nonlinear model
  predictive control}.
\newblock {\em IEE proceedings -- Control Theory and Applications},
  152(3):296--308, 2005.

\bibitem{domahidi12}
A.~Domahidi, A.~U. Zgraggen, M.~N. Zeilinger, M.~Morari, and C.~N. Jones.
\newblock Efficient interior point methods for multistage problems arising in
  receding horizon control.
\newblock In {\em Proc. of the 51\textsuperscript{st} IEEE Conf. Dec. and
  Contr. (CDC)}, pages 668--674, 2012.

\bibitem{feller15c}
C.~Feller and C.~Ebenbauer.
\newblock {Input-to-state stability properties of relaxed barrier function
  based MPC}.
\newblock In {\em Proc. of the 5\textsuperscript{th} IFAC NMPC Conf.}, pages
  302--307, Seville, Spain, 2015.

\bibitem{feller15a}
C.~Feller and C.~Ebenbauer.
\newblock Relaxed logarithmic barrier function based model predictive control
  of linear systems.
\newblock {\em arXiv:1503.03314 [math.OC]}, 2015.

\bibitem{feller15b}
C.~Feller and C.~Ebenbauer.
\newblock {Weight recentered barrier functions and smooth polytopic terminal
  set formulations for linear model predictive control}.
\newblock In {\em Proc. of the 2015 American Contr. Conf.}, pages 1647--1652,
  Chicago, IL, USA, 2015.

\bibitem{ferreau08}
H.~J. Ferreau, H.~G. Bock, and M.~Diehl.
\newblock An online active set strategy to overcome the limitations of explicit
  {MPC}.
\newblock {\em International Journal of Robust and Nonlinear Control},
  18(8):816--830, 2008.

\bibitem{giselsson14}
P.~Giselsson.
\newblock Improved fast dual gradient methods for embedded model predictive
  control.
\newblock In {\em Proc. of the 19\textsuperscript{th} IFAC World Congress},
  pages 2303--2309, 2014.

\bibitem{halanay2000}
A.~Halanay and V.~Rasvan.
\newblock {\em Stability and stable oscillations in discrete time systems}.
\newblock CRC Press, 2000.

\bibitem{hauser06}
J.~Hauser and A.~Saccon.
\newblock {A Barrier Function Method for the Optimization of Trajectory
  Functionals with Constraints}.
\newblock In {\em Proc. of the 45\textsuperscript{th} Conf. Dec. and Contr.
  (CDC)}, pages 864--869, San Diego, USA, 2006.

\bibitem{khalil02}
Hassan Khalil.
\newblock {\em Nonlinear Systems}.
\newblock Prentice Hall, 2002.

\bibitem{koegel11}
M.~K{\"o}gel and R.~Findeisen.
\newblock {Fast predictive control of linear systems combining Nesterov's
  gradient method and the method of multipliers}.
\newblock In {\em Proc. of the 50\textsuperscript{th} IEEE Conf. on Dec. and
  Contr. (CDC-ECC)}, pages 501--506, 2011.

\bibitem{mayne00}
D.~Q. Mayne, J.~B. Rawlings, C.~V. Rao, and P.~O.~M. Scokaert.
\newblock {Constrained model predictive control: Stability and optimality}.
\newblock {\em Automatica}, 36:789--814, 2000.

\bibitem{mcgovern99}
L.~K. McGovern and E.~Feron.
\newblock Closed-loop stability of systems driven by real-time, dynamic
  optimization algorithms.
\newblock In {\em Proc. of the 38\textsuperscript{th} IEEE Conf. Dec. and
  Contr. (CDC)}, volume~4, pages 3690--3696, 1999.

\bibitem{more94}
J.~J. Mor\'{e} and D.~J. Thuente.
\newblock Line search algorithms with guaranteed sufficient decrease.
\newblock {\em ACM Transactions on Mathematical Software (TOMS)},
  20(3):286--307, 1994.

\bibitem{nash94}
S.~G. Nash, R.~Polyak, and A.~Sofer.
\newblock {\em Large Scale Optimization: State of the Art}, chapter A numerical
  comparison of barrier and modified barrier methods for large-scale
  bound-constrained optimization, pages 319--338.
\newblock Kluwer Academic Publishers, 1994.

\bibitem{nocedalWright99}
J.~Nocedal and S.~J. Wright.
\newblock {\em Numerical Optimization}.
\newblock Springer, 1999.

\bibitem{patrinos12}
P.~Patrinos and A.~Bemporad.
\newblock {Simple and certifiable quadratic programming algorithms for embedded
  linear model predictive control}.
\newblock In {\em Proc. of the 4\textsuperscript{th} IFAC NMPC Conf.}, pages
  14--20, Noordwijkerhout, Netherlands, 2012.

\bibitem{rao98}
C.~V. Rao, S.~J. Wright, and J.~B. Rawlings.
\newblock Application of interior point methods to model predictive control.
\newblock {\em Journal of Optimization Theory and Applications},
  99(3):723--757, 1998.

\bibitem{richter12}
S.~Richter, C.~N. Jones, and M.~Morari.
\newblock Computational complexity certification for real-time {MPC} with input
  constraints based on the fast gradient method.
\newblock {\em IEEE Trans. Autom. Control}, 57(6):1391--1403, 2012.

\bibitem{scokaert99}
P.~O.~M. Scokaert, D.~Q. Mayne, and J.~B. Rawlings.
\newblock Suboptimal model predictive control (feasibility implies stability).
\newblock {\em IEEE Trans. Autom. Control}, 44(3):648--654, 1999.

\bibitem{tondel03}
P.~T{\o}ndel, T.~A. Johansen, and A.~Bemporad.
\newblock {An algorithm for multi-parametric quadratic programming and explicit
  {MPC} solutions}.
\newblock {\em Automatica}, 39:489--497, 2003.

\bibitem{wang10}
Y.~Wang and S.~Boyd.
\newblock {Fast Model Predictive Control Using Online Optimization}.
\newblock {\em IEEE Trans. on Control Systems Technology}, 18(2):267--278,
  2012.

\bibitem{wright1997}
S.~J. Wright.
\newblock Applying new optimization algorithms to model predictive control.
\newblock In {\em AIChE Symposium Series}, volume~93, pages 147--155, 1997.

\bibitem{zavala09}
V.~M. Zavala and L.~T. Biegler.
\newblock {The advanced-step NMPC controller: Optimality, stability and
  robustness}.
\newblock {\em Automatica}, 45(1):86--93, 2009.

\bibitem{zeilinger11}
M.~N. Zeilinger and C.~N. Jones.
\newblock {Real-Time Suboptimal Model Predictive Control Using a Combination of
  Explicit MPC and Online Optimization}.
\newblock {\em IEEE Trans. Autom. Control}, 56(7):1524--1534, 2011.

\end{thebibliography}
\end{document}